\documentclass[smallextended,final,envcountsect]{svjour3} 
\smartqed  

\usepackage{graphicx}
\usepackage{bigints}
\usepackage{todonotes}
 \presetkeys{todonotes}{color=orange!30}{}
\usepackage{color}

\usepackage{amsmath} 
\usepackage{amssymb}
\usepackage{hyperref} 
\usepackage{varioref} 

\usepackage{color}

\newcommand{\be}{\end{eqnarray*}}
\newcommand{\ee}{\end{eqnarray*}}
\newcommand{\ben}{\begin{eqnarray}}
\newcommand{\een}{\end{eqnarray}}
\newcommand{\R}{\mathbb R}
\newcommand{\N}{\mathbb N}

\newcommand{\B}{\mathbb B}
\newcommand{\C}{\mathcal C}

\newcommand{\A}{\mathcal A}

\newcommand{\T}{\mathcal T}

\newcommand{\vsm}{\vskip 0.2 truecm}
\newcommand{\ds}{\displaystyle}

\newcommand{\CC}{{\bold C}}
\newcommand{\bel}{\begin{equation}\label}
\newcommand{\eeq}{\end{equation}}
\newcommand{\w}{\omega}

\begin{document}
\title{No infimum gap and normality in optimal impulsive control under state constraints \thanks{This research is partially supported by the  INdAM-GNAMPA Project 2020 ``Extended control problems: gap, higher order conditions and   Lyapunov functions"  and by the Padua University grant SID 2018 ``Controllability, stabilizability and infimum gaps for control systems'', prot. BIRD 187147.}
}

\titlerunning{Gap and normality  in optimal impulsive control}        

\author{Giovanni Fusco        \and
        Monica Motta 
}


\institute{Giovanni Fusco\at
              Department of Mathematics "Tullio Levi-Civita", University of Padua \\ Via Trieste, 63, Padova  35121, Italy \\
              Tel.: (+39) 049 827 1261\\
              \email{fusco@math.unipd.it}           
           \and
           Monica Motta \at
              Department of Mathematics "Tullio Levi-Civita", University of Padua \\ Via Trieste, 63, Padova  35121, Italy \\
              Tel.: (+39) 049 827 1368\\
              Fax: (+39) 049 827 1499\\
              \email{motta@math.unipd.it}   
}

\date{Received: date / Accepted: date}

\maketitle

\begin{abstract}
In this paper we consider an impulsive extension of an  optimal control problem with unbounded controls, subject to endpoint and state constraints. We show that the existence of an extended-sense minimizer that is a normal extremal for a  constrained Maximum Principle ensures that there is no gap between the infima of the original problem and of its extension. Furthermore, we translate such relation into verifiable  sufficient conditions for  normality  in the form of constraint and endpoint qualifications.  Links between existence of an infimum gap and normality in impulsive control have previously been explored for problems without  state constraints. This paper establishes  such links in the presence of  state constraints and of an additional ordinary control, for locally Lipschitz continuous data.
\end{abstract}

\section*{Introduction}
In Optimal Control Theory  it is quite common practice to extend the domain of a minimum problem  to ensure the existence of the minimum  or to identify optimality conditions. In doing this, it is  of course  desirable to avoid the so-called {\em infimum gap} phenomenon, i.e. that the minimum of the extended problem is different from the minimum of the original problem.
 This is relevant not only for theoretical reasons  of `well-posedness' of the extension, but also for the actual usefulness of the extended problem in order to identify, for instance,  necessary optimality conditions  or a non degenerate Hamilton-Jacobi-Bellman  equation   for the original minimum problem. For the classical extension of a minimum problem  by convex relaxation --where the original velocity set of the trajectories is replaced by its convexification--  it has emerged that a sufficient condition to avoid the infimum gap 
is the {\em  normality} of an extended sense minimizer, namely, that {\em all} sets of multipliers verifying a Maximum Principle have cost multiplier, $\lambda$ in the following,  different from zero \cite{warga,warga1,warga2,PV1,PV2,V19}.  In  \cite{MRV}  the `normality test' has proved sufficient to guarantee the absence of an infimum gap also for the impulsive extension of an optimal control problem with unbounded controls. Very recently, in \cite{PR19,PR20} this link between normality and no-infimum-gap has been established  for the extension of an abstract  minimum control problem,
which includes both  relaxation and  impulsive extension. All the above results  in the case of the impulsive extension   concern  problems without state constraints, with $C^1$ data,  and no  ordinary controls  in the  dynamics. However , state  constraints, together with nonsmoothness of the data  and additional ordinary controls,   arise very frequently in the applications of  impulsive   optimal control  (see e.g.   \cite{MiRu,ST00,HW11,AKP15,KDPS15} and the references therein).
 
This paper provides  `normality type' sufficient conditions to avoid a gap between the  infima of the following optimization problem (P) and  the extended optimization problem ${\rm (P_e)}$ below: 
\begin{equation}
\label{Pintro}\tag{P}
\qquad\qquad\text{minimize} \,\,\, \Psi (t_1,x(t_1),t_2,x(t_2),v(t_2))
\end{equation}
over  $t_1, \, t_2\in\R$,  \,  $t_1<t_2$,  \, $(x,v, u) \in W^{1,1}([t_1,t_2];\R^{ n+1+m})$, $a\in L^1([t_1,t_2];A)$  satisfying 
  \begin{equation*}
  \begin{cases}
\dfrac{d x}{dt}(t)= f(t,x(t),a(t)) + \sum_{j=1}^{m}g_{j}(t,x(t))\,\dfrac{d u^j}{dt}(t)  \qquad \text{ a.e. } t \in [t_1,t_2],  \\[1.5ex]
\dfrac{d v}{dt}(t)=  \left|\dfrac{d u}{dt}(t)\right|  \qquad \text{ a.e. } t \in [t_1,t_2],  \\[1.5ex]
 \dfrac{d u}{dt}(t)  \in \C  \qquad\text{ a.e. } t \in [t_1,t_2], \\[1.5ex]
h_1(t, x(t)) \leq 0, \dots, h_N(t, x(t)) \leq 0  \qquad \text{ for all } t \in [t_1,t_2], \\
v(t_1)=0, \quad v(t_2)\leq K,     \quad  \left(t_1,x(t_1),t_2, x(t_2)\right) \in \T_0, 
\end{cases}
\end{equation*}
where  $K>0$ is a fixed constant (possibly equal to $+\infty$), $A\subset\R^q$ is a compact subset, $\C\subseteq\R^m$ is  a closed convex cone,   $\T_0\subseteq \R^{1+n+1+n}$ is a closed subset, and the data are locally Lipschitz continuous in  $t$,   $x$  (the precise assumptions will be given in Section \ref{S1}). Problem  (P)  is a free end-time minimization problem depending on  an ordinary control $a$   and on a control $u$  whose derivatives appear linearly in the dynamics. Furthermore, there are   time-dependent  state constraints in the form of $N$ inequalities,  endpoint constraints,  and we may have a bound $K$ on the total variation  of  $u$ --notice that  $v$ is nothing but the total variation function of $u$.  
Due to a lack of  coerciveness,  minimizers for  problem  (P) do not exist in general. Hence, adopting a by now  standard extension, we embed the original problem into the   {\it space-time} problem  ${\rm (P_e)}$ below, where  the extended state variable is  $(y^0,y,\nu):=(t,x,v)$, and the extended trajectories are $(t,x,v)$-paths  which are (reparameterized) $C^0$-limits of graphs of the original trajectories  \cite{Ris:65,War:65,BR:88,Mi:94,MR:95}:\,\footnote{As it is well-known,  a distributional approach, where  $\frac{du}{dt}$ is interpreted as a Radon measure, does not work unless  $g_i=g_i(x)$ and the Lie brackets $[g_i,g_j](x)\equiv 0$ for every $i,j=1,\dots,m$  (see e.g.  \cite{Haj85,BR:88}).}
\begin{equation}
\label{Peintro}\tag{P$_e$}
 \qquad\qquad\text{minimize } \,\,\,\Psi(y^0(0),y(0),y^0(S),y(S),\nu(S)) 
\end{equation}
over $S>0$, \, $({y^0}, y,\nu) \in W^{1,1}([0, S];\R^{1+n+1})$, \,\,\, $(\w^0,\w,\alpha)\in L^1([0, S]; \CC\times A)$  satisfying
\begin{equation*}
\begin{cases}
\dfrac{dy^0}{ds}(s) =  \w^0(s) \qquad \text{a.e. } s \in [0, S],  \\ 
\dfrac{dy}{ds}(s) =   f(y^0,y, \alpha)(s)  \w^0(s) + \sum_{j=1}^m g_j(y^0,y)(s) \w^j (s)  \quad \text{ a.e. } s \in [0, S],  \\ 
\dfrac{d\nu}{ds}(s) = |\w(s)| \qquad \text{ a.e. } s \in [0,S],  \\[1.5ex] 
h_1(y^0(s),y(s)) \leq 0, \dots, h_N(y^0(s),y(s)) \leq 0  \qquad \text{ for all } s \in [0,S], \\
\ds\nu(0)=0, \quad \nu(S)\leq K,     \quad  \left(y^0(0),y(0),y^0( S),y( S)\right) \in \T_0, 
\end{cases}
\end{equation*}
where
$
\CC:=\left\{(\w^0,\w): \ \ (\w^0,\w)\in[0,+\infty[\times\C, \ \ \w^0+|\w|=1\right\}.
$

\noindent 
 To any process $(t_1,t_2, u, a, x,v)$ of problem (P), by setting $\sigma(t):=t-t_1+v(t)$, $t\in[t_1,t_2]$, through the time-change $y^0:=\sigma^{-1}$   we  can associate a process $(S, \w^0 ,\w ,\alpha,y^0 , y ,\nu)$ of the extended problem   with  $\w^0=dy^0(s)/ds >0$  a.e.. In particular, 
 problem (P) can be identified with the restriction of  problem ${\rm (P_e)}$ to the set of processes with $\w^0>0$ a.e. (see Section \ref{S1}). Let us refer to such processes as {\em embedded strict sense processes} in the following. The extension consists therefore in considering extended sense processes with $\w^0=0$ on non-degenerate intervals, where the time $t=y^0(s)$ is constant but the extended state $y(s)$ evolves  according to  the `fast' dynamics $ dy(s)/ds =  \sum_{j=1}^m g_j(y^0(s),y(s)) \w^j (s)$.  This explains why    \eqref{Pe} is also called the  {\em impulsive problem}, although it is a conventional optimization problem  with bounded   controls. In fact, one could give an equivalent $t$-based description of this extension using bounded variation trajectories and controls   \cite{MiRu,Sa91,WZ07,KDPS15,AKP15,AR15,MS18,MS20}.

The main result of the paper, obtained  in Theorem \ref{Th_Norm} below, establishes that the existence of an extended sense minimizer for problem  ${\rm (P_e)}$ which is a normal extremal for a constrained version of the Maximum Principle, is a sufficient condition for the infimum gap avoidance.  The occurrence  of a  gap is strictly related to the  presence of endpoints and state constraints. In particular,  since the  set of trajectories corresponding to embedded strict sense processes is $C^0$-dense in the set of trajectories of the extended system, the infimum gap phenomenon can show up only when  some extended sense process  verifying the constraints is {\em isolated}, namely   cannot be approximated by  trajectories of the original system that satisfy the constraints.   From this observation,  Theorem \ref{Th_Norm}   will be derived   from a   general result  on the properties of isolated processes (see Theorem  \ref{th_isol}).  The  proof makes use of perturbation  and penalization techniques and of the Ekeland’s variational principle, in the same spirit of \cite{PV1,MRV}.  This approach is very different from that of \cite{warga,warga1,warga2,PR20}, which is based on  the construction of approximating cones to reachable sets and on set separation arguments.   We recall that  normality is not necessary to exclude the gap phenomenon: for example, it is known that without the drift $f$ in the dynamics,   gap never occurs  \cite[Lemma 4.1]{MRV} (see also   \cite{Mar00}).

The normality criterium for the absence of an infimum gap  has   some disadvantages. First of all, it requires to know a priori a minimizing extended sense  process, information that is not always available. Then,  it is necessary to verify that  all  sets of   multipliers associated to the minimizer  that meet the conditions of the Maximum Principle  have $\lambda>0$. In addition, in the presence of state constraints the normality condition may never be met, making the criterium in fact useless.  In particular, it is well known that  when the state constraint is active at the initial point of a minimizing process,  sets of degenerate multipliers  with $\lambda=0$  may always exist. Rather surprisingly, it seems that no attention has been paid to this `degeneracy question' in previous articles on the relationship between gap and normality in the presence of state constraints.

Based on the above considerations, in the second part of the paper  we first introduce  a nondegenerate version of the Maximum Principle and provide simple geometrical conditions on  endpoint  and state constraints, under which abnormal --namely,  not normal-- extremals for the original Maximum Principle  turn out to be  abnormal extremals also for the nondegenerate Maximum Principle.
 In this case, Theorem \ref{Th_Norm}   can be rephrased as follows: `normality {\it among nondegenerate multipliers} implies no infimum gap'.   This   `nondegenerate normality test'  is useful especially because in certain special cases it allows to deduce the absence of the infimum gap from easily verifiable conditions, some examples of which we  will provide. In particular, these are  constraint and endpoint  qualification conditions. 
 
Although this article is mainly focused on the infimum gap phenomenon, it also establishes  some new sufficient conditions for  normality which extend   previous conditions in \cite{MS20}.  In the literature on conventional, non-impulsive problems with state constraints,  a variety of constraint qualifications to avoid degeneracy as well as to ensure  normality are known (see e.g. \cite{FeV94,FeFoV99,FoFr15,FrTo13,LFodaP11,PV2,RV99,AA97,A00,AK16,AKP17} and  the references therein). In impulsive control, instead, some nondegenerate  Maximum Principles   have been  obtained   in \cite{AKP05,AKP15,K06,MS20}, while a Maximum Principle in normal form has only recently been introduced in \cite{MS20}.

The paper is organized as follows: in  Section \ref{S1}  we introduce precisely problems (P), ${\rm (P_e)}$ and a constrained version of the Maximum Principle for the extended problem.  Section  \ref{S2}  is devoted to we prove that an {\it  isolated}  extended sense  extremal cannot be normal and, as a corollary, we deduce that presence of an infimum gap implies abnormality of any extended sense minimizer.  In Section  \ref{S3} we provide sufficient conditions for normality, which guarantee a priori, without any knowledge of the multipliers associated with the given extended sense minimizer, the non occurrence of gap-phenomena. In   Section  \ref{S4}, we  present  some  examples to illustrate the results.

\subsection{\bf{Notations and preliminaries}}\label{sub1.1}
Given an interval $I\subseteq\R$ and a set $X \subseteq \R^k$, we write $W^{1,1}(I;X)$, $C^{0,1}(I;X)$, $C^{0,1}_{loc} (I;X)$ for the space of absolutely continuous functions, Lipschitz continuous functions, locally Lipschitz continuous functions defined on $I$ and with values in $X$,  respectively.   
 For all the classes of functions introduced so far, we will not specify domain and  codomain when the meaning is clear.   Furthermore,  we denote by $\ell(X)$, $co(X)$, Int$(X)$, $\partial X$ the Lebesgue measure, the convex hull, the interior and the boundary of $X$, respectively. As customary,  $\chi_{_X}$ is the characteristic function of $X$, namely $\chi_{_X}(x)=1$ if $x\in X$ and $\chi_{_X}(x)=0$ if $x\in\R^k\setminus X$;   $I\cdot X$ denotes the set $\{r\, x \,|\, r\in I, \, x\in X   \}$.  Given two nonempty subsets $X_1$, $X_2$ of $\R^k$, we denote by  $X_1+X_2$ the set $\{x_1+x_2 \,|\, x_1\in X_1, \, x_2\in X_2   \}$.
 Let $X\subseteq\R^{k_1+k_2}$ for some natural numbers $k_1$, $k_2$,  and write  $x=(x_1,x_2)\in\R^{k_1}\times\R^{k_2}$ for any $x\in X$. Then   $\text{{\rm proj}}_{x_i}X$ will denote the  projection of $X$ on $\R^{k_i}$,  for $i=1,2$. We denote the closed unit ball in $\R^k$ by $\B_k$, omitting the dimension when it is clear from the context. Given a closed set ${\mathcal O} \subseteq \R^k$  and a point $z \in \R^k$, we define the distance of $z$ from ${\mathcal O}$ as $d_{\mathcal O}(z) := \min_{y \in {\mathcal O}} |z-y|$.
 We set $\R_{\geq 0} := [0,+\infty[$. For any  $a,b \in \R$, we write $a \vee b:= \max \{a,b\}$.
 
 \noindent For all   $\tau_1$, $\tau_2$, $\bar \tau_1$, $\bar \tau_2\in \R $, $\tau_1<\tau_2$, $\bar\tau_1<\bar\tau_2$,  and  for any pair $(z_1, z_2)\in C^0([\tau_1,\tau_2],\R^k)\times C^0([\bar\tau_1,\bar\tau_2],\R^k)$, let us define the distance
 \begin{equation}\label{dinfty}
d_\infty\big((\tau_1,\tau_2, z_1),(\bar\tau_1, \bar\tau_2, z_2)\big) :=  
 | \tau_1 - \bar\tau_1|+  | \tau_2 - \bar\tau_2|+ \|\tilde z_1- \tilde z_2\|_{L^\infty(\R)} ,
  \end{equation} 
  where  for any  $z\in C^0([a,b],\R^k)$,  $\tilde z$  denotes its continuous constant extension to $\R$ and $\| \cdot \|_{L^{\infty}(I)}$ is the ess-sup norm on $I \subseteq \R$ interval. When the domain is clear, we will sometimes  simply write $\| \cdot \|_{L^{\infty}}$.
    
  \noindent We denote by  $NBV^{+}([0,  S];\R)$  the space of increasing, real valued functions $\mu$ on $[0, S]$ of bounded variation, vanishing at the point 0 and right continuous on $]0, S[$.   Each  $\mu\in NBV^+([0,S];\R)$  defines a Borel measure on $[0, S]$,   still denoted by $\mu$,  its total variation function is  indicated by  $\| \mu \|_{TV}$ or equivalently by $\mu([0, S])$, and   its support  by spt$\{\mu\}$.  
 
\vsm 
Some standard constructs from nonsmooth analysis are employed in this paper. For background material we refer the reader for instance to \cite{OptV}. A set 
$K \subseteq \R^k$ is a {\em cone} if $\alpha k \in K$ for any $\alpha >0$,  whenever $k \in K$. Take a closed set $D \subseteq \R^k$ and a point $\bar x \in D$, the \textit{proximal normal cone} $N^P_D(\bar x)$ of $D$ at $\bar x$ is defined as
\[
N^P_D(\bar x) := \left\{  \eta\in\R^k \text{ : } \exists M>0 \,\,\text{ such that }\,\, \eta \cdot (x-\bar x) \leq M |x-\bar x|^2 \,\, \forall x \in D     \right\}.
\]
The \textit{limiting normal cone} $N_D(\bar x)$ of $D$ at $\bar x$ is given by 
\[ 
N_D(\bar x) := \left\{  \eta\in\R^k \text{ : } \exists x_i \stackrel{D}{\to} \bar x,\, \eta_i \to \eta \,\,\text{ s.t. }\,\,  \eta_i \in N^P_D(x_i) \text{ for each } i\in \N     \right\},
\]
 in which the notation $x_{i} \stackrel{D}{\longrightarrow}\bar{x}$ is used to indicate that  all points in the converging sequence  $(x_i)_{i\in\N}$  lay in $D$.   In general,  $N^P_D(\bar x) \subseteq N_D(\bar x)$.  
Take a lower semicontinuous function  $G:\R^k \to \R$  and a point $\bar x \in \R^k$, the \textit{limiting subdifferential} of $G$ at $\bar x$ is 
\[
\partial G(\bar x) := \left\{     \xi \in \R^k \text{ : } (\xi,-1) \in N_{epi(G)}(\bar x, G(\bar x)) \right\},
\]
where $epi(G)$ is the \textit{epigraph} of $G$. If  $G:\R^k\times\R^{h} \to \R$ is a lower semicontinuous function and $(\bar x,\bar y) \in \R^k\times\R^{h}$, we  write $\partial_x G(\bar x,\bar y)$, $\partial_y G(\bar x,\bar y)$ to denote the {\em partial limiting subdifferential of $G$ at  $(\bar x,\bar y)$ w.r.t. $x$, $y$}, respectively. 
 Given $G \in C^{0,1}_{loc} (\R^k;\R)$ and $\bar x \in \R^k$, the \textit{reachable hybrid subdifferential} of $G$ at $\bar x$ is
$$
  \partial^{*>} G(\bar x):= \left\{ \xi \in \R^k\,: \, \exists\; 
x_{i} \rightarrow \bar{x} \mbox{ s. t.  } \ G(x_i)>0 \ \text{$\forall i$ and} \ \ \nabla G(x_i)\to\xi\right\},
$$
while the \textit{reachable gradient} of $G$ at $\bar x$ is 
\[
\partial^{*} G(\bar x) := \, \left\{    \xi \in \R^k \text{ : } \exists (x_i)_i  \subset \text{diff}(G) \setminus \{ \bar x\} \text{ s.t. }  x_i \to \bar x  \text{ and }   \nabla G(x_i) \to \xi       \right\}
\] 
where diff($G$) denotes the set of differentiability points of $G$ and $\nabla$ is  the usual gradient operator.
We define the \textit{hybrid subdifferential} as $\partial^{>} G(\bar x):=$co$\,\partial^{*>} G(\bar x)$. The set $\partial^*G(\bar x)$ is nonempty, closed, in general non convex, and its convex hull coincides with the Clarke subdifferential $\partial^{C} G(\bar x)$, that is $\partial^{C} G(\bar x)=$co$\,\partial^{*} G(\bar x)$. Finally, when $G$ is locally Lipschitz continuous, $\partial^{C} G(\bar x)=$co$\,\partial  G(\bar x)$.




\section{Optimal control problems and a Maximum Principle}\label{S1}
In this section we introduce rigorously the constrained optimization problem over $W^{1,1}$-controls $u$ and its embedding in an extended, or impulsive, problem. Furthermore, we  state a Maximum Principle for the extended problem. For simplicity, we will establish all the results for a single state constraint, explaining from time to time with remarks how to adapt these results to the case with $N$ constraints. 
\vsm
 Throughout the  paper we shall consider the following hypotheses. 
{\em
\begin{itemize} 
\item[] {\bf (H0)}
The control set $\C\subseteq\R^m$ is a convex  cone, the set of ordinary controls $A \subset \R^q$ is compact, and the endpoint constraint set $ \T_0\subseteq \R^{1+n+1+n}$ is closed. 
\item[] {\bf (H1)}The drift function $f \in C(\R^{1+n} \times A; \R^n)$ and, for every $a\in A$, $f(\cdot, a)  \in C_{loc}^{0,1}(\R^{1+n}; \R^n)$, uniformly w.r.t. $A$;  $g_j \in C_{loc}^{0,1}(\R^{1+n}; \R^n)$ for any $j=1,\dots,m$;  $h \in C_{loc}^{0,1}(\R^{1+n}; \R)$.
\item[] {\bf (H2)}
The cost function $\Psi \in C^{0,1}((\bar t_1,\bar x_1,\bar t_2,\bar x_2,\bar v_2)+\delta\B;\R)$ for some $\delta >0$, where $(\bar t_1,\bar x_1,\bar t_2,\bar x_2,\bar v_2)$ denotes the endpoints  of the  optimal trajectory  that we will consider in all our results; moreover,   for every $(t_1,x_1,t_2,x_2)$, the map $v_2\mapsto \Psi(t_1,x_1,t_2,x_2,v_2)$ is monotone non-decreasing.
\end{itemize}}

\subsection{\bf{The original optimal control problem}}\label{sub2.1}
We  set  $\T:= \T_0  \times ]-\infty,K]$ and define the {\em set ${\mathcal U}$ of  strict sense controls} as
$$
{\mathcal U}:=\left\{\begin{array}{l}
(t_1,t_2,u,a):  t_1,t_2\in\R, \ t_1<t_2, \  (du/dt,a)\in  L^1([t_1,t_2];\R^{m}\times A),   \\[1.5ex]
\quad\qquad\qquad\qquad\qquad  \ d u(t)/dt \in \C   \ \text{ a.e.  } \ t\in [t_1,t_2]
\end{array}\right\}. 
$$
\begin{definition}[Strict sense processes]  Let $(t_1,t_2,u,a)\in{\mathcal U}$ be a strict sense control, we call  $(t_1,t_2,  u,a, x ,v  )$ a {\em strict sense process} if  the pair $(x,v)\in  W^{1,1}([t_1,t_2];\R^{n+1})$ verifies
\[
\begin{cases}
\displaystyle\dfrac{d x}{dt}(t)=f(t,x(t),a(t))+\sum_{j=1}^m g_j(t,x(t)) \dfrac{d u^j}{dt}(t)   \\[1.5ex]
\dfrac{d v}{dt}(t)= \left| \dfrac{d u}{dt}(t) \right| 
\end{cases}  \qquad \text{a.e. } t\in [t_1,t_2].
 \]
Furthermore, we say that $(t_1,t_2,  u, a, x ,v )$ is  \textit{feasible}  if   $v(t_1)=0$, $h(t,x(t)) \leq 0$ for each $t\in [t_1,t_2]$ and $\big(t_1,x(t_1),t_2, x(t_2),v(t_2)\big) \in \T$.
\end{definition}
 The original optimal control problem  is defined as
\begin{equation}
\label{P}\tag{P}
\left\{\begin{array}{l}
\ds\text{minimize}  \,\,\, \Psi(t_1,x(t_1),t_2,x(t_2),v(t_2)) \\[1.5ex] 
\text{over the set of feasible strict sense processes $(t_1,t_2 , u,a, x ,v )$.}
\end{array}\right.
\end{equation}
We  consider the following concept of local minimizer.
 \begin{definition}  We call a feasible  strict sense process $(\bar t_1,\bar t_2 , \bar u,\bar a, \bar x, \bar v )$   a {\em local strict sense  minimizer} of  $\eqref{P}$ if 
 there exists $\delta >0$ such that
\begin{equation}\label{min1}
\Psi(\bar t_1,\bar x(\bar t_1),\bar t_2,\bar x(\bar t_2),\bar v(\bar t_2))\,\leq\,\Psi(t_1,x(t_1),t_2,x(t_2),v(t_2))
\end{equation}
for any feasible $(t_1,t_2 , u,a, x ,v )$ 
verifying
$
d_\infty\Big((\bar t_1,\bar t_2, \bar x, \bar v) , (t_1,t_2, x ,v) \Big)<\delta,
$
where $d_\infty$ is the distance defined in \eqref{dinfty}.
 If  relation \eqref{min1} is satisfied for all feasible strict sense processes,  
 we say that $(\bar t_1,\bar t_2,  \bar u,\bar a, \bar x, \bar v   )$ is a {\em (global) strict sense minimizer}.
\end{definition}

\begin{remark}\label{Rext1} Arguing similarly to  \cite{MS20}, we could consider a more general cost  of the form 
$$
\Psi(t_1,x(t_1),t_2,x(t_2),v(t_2))+\int_{t_1}^{t_2} \ell_0(t,x(t),a(t))+\ell_1(t,x(t),|\dot u(t)|) \,dt,
$$
where    $\ell_0$, $\ell_1$  are nonnegative  and the  {\em extended Lagrangian} $L$,  defined by
 $$L(t,x,\w^0,r,a):= \ell_0(t,x,a)\w^0+ 
\lim_{\rho\to \w^0}\ell_1\left(t,x,\rho^{-1}r\right)\rho, $$
verifies $ L \in C( \R\times\R^n\times\R_{\ge0}\times\R_{\ge0}\times A; \R_{\ge0})$ and, for every $a\in A$,   $L(\cdot, a)  \in C_{loc}^{0,1}( \R\times\R^n\times\R_{\ge0}\times\R_{\ge0}; \R_{\ge0})$, uniformly w.r.t. $A$. 
The results of this article can also be applied to the case where dynamics, cost, and constraints depend on the variable $u$. In fact, it is sufficient to add to the control system in $(P)$ the equations    $d x_{n+1}(t)/dt=d u_1(t)/dt, \dots, dx_{n+m}(t)/dt=d u_m(t)/dt$. 
\end{remark}

\begin{remark} \label{constraint} 
As is not difficult to see, given a closed,  Hausdorff-Lipschitz continuous multifunction $X: [t_1, t_2] \rightsquigarrow \R^n$,    the function $h(t,x) := d_{X(t)} (x)$  belongs to $C^{0,1}(\R^{1+n};\R)$. Therefore,  we could allow implicit time-dependent state constraints of the form $x(t)  \in X(t)$   for all $t\in [t_1,t_2]$,   
since one clearly has that   $x(t)\in X(t)$  if and only if $h(t, x(t))\le0$ on  $ [t_1,t_2]$.

\end{remark}

\subsection{\bf{The  extended optimal control problem}}\label{sub2.2}
We set 
\begin{equation}\label{TCC}
 \CC:= \left\{ (\omega^0,\omega) \in \R_{\geq 0}\times \C \text{ : } \omega^0 + |\omega|=1 \right\} 
 \end{equation}
and introduce the   {\em set of extended sense controls,} defined as follows:
$$
\mathcal{W}:=\bigcup_{S>0} \{S\}\times  L^1([0,S]; \CC\times A). 
 $$ 
 
 \begin{definition}[Extended sense processes]  For any extended sense control   $(S, \w^0,\w,\alpha)\in \mathcal{W},$  we refer to $(S, \w^0 ,\w ,\alpha,y^0 , y ,\nu)$ as an   \textit{extended sense  process}  if   $(y^0,y,\nu)\in  W^{1,1}([0,S];\R^{1+n+1}$ verifies
 \begin{equation}
 \label{extended}
 \left\{
\begin{split}
\frac{d{y^0}}{ds} (s) &= \w^0(s)  \\
\frac{dy}{ds} (s) & = f(y^0, y,\alpha)(s)\w^0(s)+ \sum_{i=1}^{m}g_{i}(y^0, y)(s)\w^i(s)\\
\frac{d\nu}{ds} (s) & = |\w(s)|
\end{split}
\right. \quad {\rm a.e.}\, s\in [0,S]. \end{equation}
We say that $(S , \w^0 ,\w ,\alpha,y^0 , y ,\nu)$ is {\em  feasible} if $\nu(0)=0$, $h(y^0(s),y(s)) \leq 0$ for each $s\in [0,S]$ and  $\left(y^0(0),y(0),y^0(S), y(S),\nu(S)\right) \in \T$.
\end{definition}

 The set of strict sense processes, say $\Sigma$,  can be embedded into the set of extended sense processes, $\Sigma_e$, through the following map ${\mathcal I}: \Sigma\to\Sigma_e$, defined as 
\bel{eqex}
 {\mathcal I}(t_1,t_2, u,a, x ,v):= (S,\w^0,\w,\alpha,y^0,y,\nu),
\eeq
 where,  setting $\sigma(t):=t-t_1+v(t)$, $S:=\sigma(t_2)$, and $y^0:=\sigma^{-1}:[0,S]\to[t_1,t_2]$,
we associate to any strict sense process  $(t_1,t_2, u,a, x ,v)$  the extended sense process 
\bel{eqex}
(S,\w^0,\w,\alpha,y^0,y,\nu):=\left(S, \frac{dy^0}{ds},  \left(\frac{du}{dt} \circ y_{0}\right) \cdot \frac{dy^0}{ds},  a\circ y^{0}, y^{0},x\circ y^{0}, v\circ y^{0}\right),\, \footnote{Since  every   $L^1$-equivalence class contains  Borel measurable representatives,   we are tacitly assuming that all $L^1$-maps we are considering are Borel measurable.
}
\eeq
where clearly $\w^0>0$ a.e..  Conversely, if $(S,\w^0,\w,\alpha,y^0,y,\nu)$ is an extended sense process with   $\w^0>0$ a.e.,     the absolutely continuous, increasing and surjective inverse  $\sigma:[t_1,t_2]\to [0,S]$ of  $y^0$,  allows us  to define the strict sense process
\bel{eqex1}
(t_1,t_2, u,a,x,v):=\left(t_1,t_2, \int_{\sigma(t_1)}^{\sigma(t)}\w(s)\,ds, \alpha\circ\sigma, y\circ\sigma,\nu\circ\sigma\right). 
\eeq 
Therefore, ${\mathcal I}$ is injective, \footnote{Of course, up to translations of $u$.} ${\mathcal I}(\Sigma)=\{(S,\w^0,\w,\alpha,y^0,y,\nu)\in \Sigma_e: \  \w^0>0 \ \text{a.e.}\}$,  and   the extension consists in considering also  extended sense processes with $w^0$ possibly zero on some non-degenerate intervals. As anticipated in the Introduction, we will sometimes refer to the processes in ${\mathcal I}(\Sigma)$ as {\em embedded strict sense processes}. 
 
\vsm
We define the {\it extended problem} as
\begin{equation}
\label{Pe}\tag{P$_e$}
\left\{\begin{array}{l}
\text{minimize } \,\,\, \Psi(y^0(0),y(0),y^0(S), y(S),\nu(S))  \\[1.5ex]
\text{over  feasible estended-sense processes 
$(S , \w^0 ,\w ,\alpha,y^0 , y ,\nu)$.}
\end{array}\right.
\end{equation}

 \begin{definition}\label{emin} A feasible extended sense process   $(\bar S, \bar \w^0,\bar \w,\bar\alpha,  \bar y^0,\bar  y,\bar \nu)$   is said to be a  {\em local minimizer for the extended problem \eqref{Pe}} if 
  there exists $\delta >0$ such that
\begin{equation}
\label{min}
  \Psi(\bar y^0(0),\bar y(0),\bar y^0(S), \bar y(S),\bar \nu(S))) \leq  
   \Psi(y^0(0),y(0),y^0(S), y(S),\nu(S))
\end{equation}
 for all feasible extended sense processes $(S, \w^0 ,\w ,\alpha,y^0 , y ,\nu)$ that satisfy \linebreak
$
d_\infty\Big((y^0(0),y^0(S), y ,\nu) , (\bar y^0(0),\bar y^0(\bar S), \bar y ,\bar \nu)\Big)$ $<\delta,
$
where $d_\infty$ is as  in \eqref{dinfty}.  If \eqref{min} is satisfied for all feasible extended sense processes, 
we call
 $(\bar S, \bar \w^0,\bar \w,\bar\alpha,  \bar y^0,\bar  y,\bar \nu)$   a {\em (global) extended sense  minimizer}. 
\end{definition}

\begin{remark}\label{Lext}  The notion of extended sense local minimizer is the natural extension of the definition of strict sense local minimizer. Indeed, in view of  \cite[Prop. 2.7]{AMR19},   $(\bar T, \bar u, \bar a, \bar x,\bar v)$ is a strict sense  local minimizer for \eqref{P}  if and only if   ${\mathcal I}(\bar T, \bar u, \bar a, \bar x,\bar v)$ is an  extended sense   local minimizer  for \eqref{Pe}  among  the  feasible  embedded strict sense processes.  \end{remark}

%
%



\subsection{\bf{A Maximum Principle for the extended problem}}
Consider  the  {\it unmaximized Hamiltonian} $H$, defined  by
$$
\begin{array}{l}
H(t,x,p_0,p,\pi,\w^0,\w,a)   := p_0\w^0 + p\cdot\Big(f(t,x,a) \w^0 +  \sum_{i=1}^{m}  g_{i}(t,x) \w^i\Big) + \pi | \w|,
\end{array}
$$ 
for all $(t,x,p_0,p,\pi,\w^0,\w,a)\in \R^{1+n+1+n+1 }\times \CC\times A$.

\begin{theorem}[PMP] \label{Th_pmp}
 Assume  {\bf (H0)}-{\bf (H2)}.
Let $(\bar S, \bar\w^0,\bar\w,\bar\alpha, \bar{y}^0,\bar y,\bar\nu)$ be an extended sense local  minimizer for \eqref{Pe}. Then  there exist a path $(p_0,p) \in W^{1,1}([0, \bar S];\R\times\R^n)$, $\lambda \geq 0$, $\pi \leq 0$, $\mu \in NBV^{+}([0, \bar S];\R)$, $(m_0,m) : [0, \bar S] \to \R^{1+n}$ Borel measurable and $\mu$-integrable functions, verifying the following conditions:
\begin{itemize}
\item[{\rm (i)}] {\sc (non-triviality)} 
\begin{equation}\label{fe1}
\| p_0 \|_{L^{\infty}} +\| p \|_{L^{\infty}}+ \| \mu \|_{TV} + \lambda \neq 0; 
\end{equation}
\item[{\rm (ii)}]  {\sc (adjoint equation)} for a.e. $s \in [0, \bar S]$,
\[\left(-\frac{dp_0}{ds}(s), -\frac{dp}{ds}(s)\right) \in co \ \partial_{t,x} \ H\Big(\bar y^0(s), \bar y(s), q_0(s) , q(s),\pi,  \bar \w^0(s),\bar \w(s),\bar\alpha(s)\Big), \]
\item[{\rm (iii)}] {\sc  (tranversality)} 
\begin{equation} \label{trans_cond}
\begin{split}
&\left(p_0(0), p(0), -q_0(\bar S), -q(\bar S) , -\pi \right) \in \lambda \partial \Psi\left(\bar{y}^0(0),\bar{y}(0),\bar{y}^0(\bar S), \bar y(\bar S), \bar\nu(\bar S)\right) \\ 
&\textcolor{white}{d}\qquad\qquad\qquad\qquad\qquad\qquad \qquad +N_{\T}\left(\bar{y}^0(0),\bar{y}(0),\bar{y}^0(\bar S), \bar y(\bar S),   \bar\nu(\bar S)\right);
\end{split}
\end{equation}
\item[{\rm (iv)}] {\sc (maximization and vanishing of the Hamiltonian)} for a.e. $s \in [0, \bar S]$,
\begin{equation} \label{maxham}
\begin{array}{l}
H\Big(\bar y^0(s), \bar y(s), q_0(s) , q(s),\pi,  \bar \w^0(s),\bar \w(s),\bar\alpha(s)\Big) \\
\ds \qquad\qquad\qquad=  \max_{(\w^0,\w,a)\in \CC\times A}  H\Big(\bar y^0(s), \bar y(s), q_0(s) , q(s),\pi,  \w^0,\w,a\Big)=0;  
\end{array}
\end{equation}
\item[{\rm (v)}] 
$(m_0,m)(s) \in \partial_{t,x}^{>}\, h\left(\bar{y}^0(s), \bar y(s)\right)$ \qquad $\mu$-a.e. $s \in [0, \bar S]$;

\item[{\rm (vi)}] $spt(\mu) \subseteq \{ s\in[0,\bar S] \text{ : } h\left( \bar{y}^0(s), \bar y(s)\right) = 0 \}$,
\end{itemize}
where
\[
(q_0,q)(s) := 
\begin{cases}
\ds (p_0,p)(s) +  \int_{[0,s[} (m_0,m)(\tau) \mu(d\tau)\, \qquad \qquad \qquad \,\, s\in[0, \bar S[, \\
\ds (p_0,p)(\bar S) + \int_{[0, \bar S]} (m_0,m)(\tau) \mu(d\tau)\, \qquad \qquad\qquad s=\bar S.
\end{cases}
\]
Furthermore:
\begin{itemize}
\item[{\rm (vii)}] if $\lambda \partial_v\Psi\left(\bar{y}^0(0),\bar{y}(0),\bar{y}^0(\bar S), \bar y(\bar S), \bar\nu(\bar S)\right) = 0$ and $\bar\nu(\bar S) < K$, then $\pi=0$;

\item[{\rm (viii)}] if $\bar y^0(0) < \bar y^0(\bar S)$, (i) can be strengthened to $\| p \|_{L^{\infty}}+ \| \mu \|_{TV} + \lambda \neq 0$.
\end{itemize}
\end{theorem}
\noindent {\it Proof: }\,  The extended problem \eqref{Pe} is a conventional  optimization problem in the state-space $(y^0,y,\nu)\in\R^{1+n+1}$, with   endpoint constraint $\T$ and   state constraint $h(y^0,y)\le0$, to which  standard  `free end-time'  versions of the constrained Maximum Principle  are applicable.   In particular, the current result  can be  deduced from \cite[Theorem 9.3.1]{OptV} by means of usual reparameterization techniques (see e.g. \cite[Theorem 8.7.1]{OptV}). Actually, by these arguments   it follows the existence of a  further multiplier $r \in W^{1,1}([0,\bar S]; \R)$ such that
\[
\begin{aligned}
& \frac{dr}{ds}(s) = -H_s\Big(\bar y^0(s), \bar y(s), q_0(s) , q(s),\pi,  \bar \w^0(s),\bar \w(s),\bar\alpha(s)\Big) \quad &\text{a.e. }s\in[0,\bar S],  \\
& r(\bar S) \in -\lambda\Psi_s \left(\bar{y}^0(0),\bar{y}(0),\bar{y}^0(\bar S), \bar y(\bar S), \bar\nu(\bar S)\right),& \\
& \max_{(\w^0,\w,a)\in \CC\times A}  H\Big(\bar y^0(s), \bar y(s), q_0(s) , q(s),\pi,  \w^0,\w,a\Big)=r(s) &\text{a.e. }s\in[0,\bar S],
\end{aligned}
\]
where $H_s$, $\Psi_s$ denote the partial derivatives of $H$ and $\Psi$ w.r.t. $s$, respectively.
However, since the vector fields $f$ and $(g_j)_{j=1,\dots,m}$, the cost function $\Psi$, and the constraints do not depend  explicitly on the pseudo-time $s$, this yields the constancy of the Hamiltonian with constant equal to 0 in (iv).  
Finally,  the strengthened  non-triviality condition \eqref{fe1}, which does not involve the multiplier $\pi$ associated to  $\nu$,  and  the refinements  (vii),  (viii),   can be proved  as in \cite[Theorem 3.1]{MRV}. 
 \qed
 \vsm
 
 \begin{remark}[Multiple state constraints]\label{RMSC1}
As observed in \cite[Section 9]{OptV}, when in problem \eqref{Pe} the single state constraint  is replaced with a collection of constraints $h_i(t,x)\le 0$, \, $h_i\in C^{0,1}_{loc}(\R^{1+n};\R)$ for $i=1,\dots, N$,  from Theorem \ref{Th_pmp} one  can deduce  the following corollary:  given an extended sense local  minimizer $(\bar S, \bar\w^0,\bar\w,\bar\alpha, \bar{y}^0,\bar y,\bar\nu)$, there exist $(p_0,p) \in W^{1,1}$, $\lambda \geq 0$, $\pi \leq 0$, $\mu_i \in NBV^+([0, \bar S];\R)$ for $i= 1,\dots, N$ and Borel measurable and $\mu_i$-integrable functions $(m_{0_i},m_i)$,  such that
$
(m_{0_i},m_i)(s) \in \partial^>_{t,x}h_i\left(\bar{y}^0(s), \bar y(s)\right) \ \ \mu_i\text{-a.e. }s \in [0, \bar S]$, \,
$spt(\mu_i) \subseteq \{ s\in[0,\bar S] \text{ : } h_i\left( \bar{y}^0(s), \bar y(s)\right) = 0 \},$
and conditions (i)--(iv), (vii) and (viii) of Theorem \ref{Th_pmp} are met with $(q_0, q)$ and $\mu$  verifying
\begin{equation}\label{modifica_state_constraints}
\begin{split}
&(q_0,q)(s) := 
\begin{cases}
 (p_0,p)(s) + \int_{[0,s[} \sum_{i=1}^N (m_{0_i},m)(\tau) \mu_i(d\tau)  \quad\,\,\,\, \text{if }s<\bar S \\
(p_0,p)(\bar S) + \int_{[0,\bar S]} \sum_{i=1}^N (m_{0_i},m)(\tau) \mu_i(d\tau) \,\,\quad \text{if }s=\bar S,
\end{cases} \\
& \mu(A) = \mu_1(A)+ \dots + \mu_N(A) \qquad \text{for all Borel subsets } A\subseteq [0,\bar S] .
\end{split}
\end{equation}
\end{remark}

\begin{definition}[Normal and abnormal extremal]
We say that a   feasible extended sense
 process $(\bar S,\bar\w^0,\bar\w,\bar\alpha,\bar{y}^0,\bar y,\bar\nu)$  is   an (extended sense)  {\em  extremal}  if there exists a set of multipliers $(p_0,p,\pi, \lambda, \mu)$ and functions $m_0$ and $m$ which meet the conditions  of Theorem \ref{Th_pmp}.   We will call an extremal  \textit{normal} if all possible choices of multipliers as above have $\lambda >0$, and  \textit{abnormal}  when  there exists at least one set of such multipliers with $\lambda =0$. 
 \end{definition}


\section{Infimum gap and abnormality}\label{S2}
Write  $J(t_1,t_2, u, a, x , v):=\Psi(t_1, x(t_1), t_2, x(t_2) , v(t_2))$ for the cost of a strict sense process $(t_1,t_2, u, a, x , v)$ in problem \eqref{P},  
 and  $J_{e}(S, \w^0,\w,\alpha, y^0 , y , \nu):= \Psi(y^0(0), y(0),y^0(S), y(S) ,\nu(S))$ for the cost of an extended sense process \linebreak $(S, \w^0,\w,\alpha, y^0 , y , \nu)$ in problem \eqref{Pe}.
We also write $\Sigma^f\subseteq\Sigma$ and $\Sigma_e^f\subseteq\Sigma_e$ for the subset of feasible strict sense processes  and for  the subset of feasible extended sense processes, respectively.  
\begin{definition}[Infimum gap] \label{gapdef}
We shall say that   
 there is  {\em infimum gap} if 
$$
\ds\inf_{\Sigma_e^f}  \,J_{e}(S, \w^0,\w,\alpha, y^0 , y , \nu) < \inf_{\Sigma^f} \, J(t_1,t_2, u, a, x , v).
$$
Furthermore, if $(\bar S, \bar\w^0 ,\bar\w,\bar\alpha , \bar y^0 ,\bar y , \bar\nu )$ is an  extended sense  local  minimizer, we shall say that 
  there is  {\em  local infimum   gap at $(\bar S, \bar\w^0 ,\bar\w,\bar\alpha , \bar y^0 ,\bar y , \bar\nu )$}    if, for some $\delta>0$, 
$$
 \ds J_{e}(\bar S, \bar\w^0 ,\bar\w,\bar\alpha , \bar y^0 ,\bar y , \bar\nu )  <\,  \inf_{B^\delta (\bar S, \bar\w^0 ,\bar\w,\bar\alpha , \bar y^0 ,\bar y , \bar\nu )\cap\Sigma^f }  \, J(t_1,t_2, u, a, x , v),  
 $$
where we have set 
$$
\begin{array}{l} 
B^\delta (\bar S, \bar\w^0 ,\bar\w,\bar\alpha , \bar y^0 ,\bar y , \bar\nu )   := 
\Big\{(t_1,t_2, u, a, x , v) \in \Sigma: \ (S,\w^0,\w,\alpha, y^0 , y , \nu)  \\ 
 \quad :={\cal I} (t_1,t_2, u, a, x , v) \hbox{ and } d_\infty\Big((y^0(0),y^0(S), y ,\nu) , (\bar y^0(0),\bar y^0(\bar S), \bar y ,\bar \nu)\Big)<\delta\Big\}.
\end{array}
$$ 
 \end{definition}      
To prove that, in the presence of a gap,   extended sense  local minimizers for problem \eqref{Pe} are abnormal extremals, it is convenient to rephrase  Definition \ref{gapdef}  only  in terms of extended sense processes.  Precisely, using the above notation, by the properties of the map ${\mathcal I}$ (see \eqref{eqex})  it follows that 
$$
\ds \inf_{ \Sigma^f} \, J(t_1,t_2, u, a, x , v)=\inf_{{\mathcal I}(\Sigma^f)} \,J_{e}(S, \w^0,\w,\alpha, y^0 , y , \nu),
$$
\[
\begin{split}
&\inf_{B^\delta(\bar S, \bar\w^0 ,\bar\w,\bar\alpha , \bar y^0 ,\bar y , \bar\nu )\cap\Sigma^f }  \, J(t_1,t_2, u, a, x , v) \\
&\textcolor{white}{kjhgfldhlkjhglkdjfh}=\inf_{B^\delta_e (\bar S, \bar\w^0 ,\bar\w,\bar\alpha , \bar y^0 ,\bar y , \bar\nu )\cap{\mathcal I}(\Sigma^f)}  \,J_{e}(S, \w^0,\w,\alpha, y^0 , y , \nu),
\end{split}
\]
where
$$
\begin{array}{l} 
B^\delta_e (\bar S, \bar\w^0 ,\bar\w,\bar\alpha , \bar y^0 ,\bar y , \bar\nu )   := 
\Big\{(S,\w^0,\w,\alpha, y^0 , y , \nu)\in   \Sigma_e: \\
 \quad\qquad\qquad \qquad\qquad\qquad   d_\infty\Big((y^0(0),y^0(S), y ,\nu) , (\bar y^0(0),\bar y^0(\bar S), \bar y ,\bar \nu)\Big)<\delta\Big\}.
   \end{array}
$$ 
Even if the set of embedded strict sense processes is dense into the set of extended sense processes with respect to the distance $d_\infty$, the  infimum gap can actually occur,   since all embedded strict sense processes close to a given feasible extended sense process might violate either the endpoint constraints or the state constraint. This leads us to the following definition:

\begin{definition}[Isolated feasible extended sense process] \label{isolated}
A feasible extended sense process $(S,\w^0,\w,\alpha, y^0 , y , \nu)$ is called \textit{isolated} if,  for some $\delta>0$, one has
$$
B^\delta_e (S,\w^0,\w,\alpha, y^0 , y , \nu) \cap {\mathcal I}(\Sigma^f)=\emptyset.
$$
\end{definition}

The following result relates isolated feasible extended sense processes and  infimum gap. 
\begin{proposition} \label{iso} Assume  {\bf (H0)}-{\bf (H1)}.  
Let $(\bar S, \bar\w^0 ,\bar\w,\bar\alpha , \bar y^0 ,\bar y , \bar\nu )$ be an  extended sense  minimizer [resp., local minimizer]  for the extended problem \eqref{Pe} and assume $\Psi\in C^0((\bar y^0(0),\bar y(0), \bar y^0(\bar S) ,\bar y(\bar S) , \bar\nu(\bar S))+\delta\B;\R)$ for some $\delta>0$.  If
 there is  infimum gap  [resp., local infimum gap at  $(\bar S, \bar\w^0 ,\bar\w,\bar\alpha , \bar y^0 ,\bar y , \bar\nu )$], then  $(\bar S, \bar\w^0 ,\bar\w,\bar\alpha , \bar y^0 ,\bar y , \bar\nu )$ is an  isolated feasible extended sense process.
\end{proposition}

\noindent {\it Proof: }\, 
 Suppose by contradiction that  $(\bar S, \bar\w^0 ,\bar\w,\bar\alpha , \bar y^0 ,\bar y , \bar\nu )$ is not isolated. Then we can take a sequence $\delta_j  \downarrow 0$ and, for each $j \in \N$, there exists \linebreak $\left(S_j, \w^0_j, \w_j ,\alpha_j , y^0_j,y_j,\nu_j\right) \in B^{\delta_j}_e (\bar S, \bar\w^0 ,\bar\w,\bar\alpha , \bar y^0 ,\bar y , \bar\nu ) \cap {\mathcal I}(\Sigma^f)$. By the definition of $d_{\infty}$ and  the continuity of the cost function $\Psi$, this implies that  no infimum gap may occur.
\qed
\vsm

In the following theorem we establish the main  result of this section: 
\begin{theorem}\label{th_gap} Assume  {\bf (H0)}-{\bf (H1)}.  
Let $(\bar S, \bar\w^0 ,\bar\w,\bar\alpha , \bar y^0 ,\bar y , \bar\nu )$ be an  extended sense  minimizer [resp., local minimizer]  for the extended problem \eqref{Pe} and assume $\Psi\in C^0((\bar y^0(0),\bar y(0), \bar y^0(\bar S) ,\bar y(\bar S) , \bar\nu(\bar S))+\delta\B;\R)$ for some $\delta>0$. 
If  there is  infimum gap  [resp., local infimum gap at  $(\bar S, \bar\w^0 ,\bar\w,\bar\alpha , \bar y^0 ,\bar y , \bar\nu )$], then $(\bar S, \bar\w^0 ,\bar\w,\bar\alpha , \bar y^0 ,\bar y , \bar\nu )$ is an abnormal extremal.
\end{theorem}

Thanks to Proposition \ref{iso}, Theorem \ref{th_gap}  is a  straightforward  consequence of the following result, which  extends  \cite[Th. 4.4]{MRV} to the case with state constraints, an additional ordinary control in the drift,  and nonsmooth data. 
 
\begin{theorem}\label{th_isol} Assume  {\bf (H0)}-{\bf (H1)}.
If $(\bar S, \bar\w^0 ,\bar\w,\bar\alpha , \bar y^0 ,\bar y , \bar\nu )$ is an  isolated feasible extended sense process, then it is an abnormal extremal.
\end{theorem}
\noindent {\it Proof: }\,   
Since  the proof  involves  only space-time trajectories $(y^0,y)$  which are  close to  the reference space-time trajectory $(\bar y^0 , \bar y )$ and the controls assume values in a compact set,  using standard  truncation and mollification arguments   we can assume  that  there exists some $L>0$  such that the functions $f$, $g_1,\dots,g_m$, and $h$   are $L$-Lipschitz continuous and  bounded by $L$. 
The proof is divided into several steps in which successive sequences of optimization problems are introduced that have as eligible controls only embedded strict controls, and costs that measure how much a process violates the constraints. Using the Ekeland Principle, minimizers are then built for these problems, which converge to the initial isolated process. Furthermore,  applying the PMP to these approximate problems with reference to the above mentioned minimizers, we obtain in the limit a set of multipliers with $\lambda=0$ for problem \eqref{Pe}, with reference to the isolated process $(\bar S, \bar\w^0 ,\bar\w,\bar\alpha , \bar y^0 ,\bar y , \bar\nu )$.

\vsm
{\it Step 1.}  Define the function $\Phi:\R^{1+n+1+n+1}\to\R$,  given by
\[
\Phi \left(t_1,x_1,t_2 ,x_2,v_2\right) := d_{\T_0}(t_1,x_1,t_2,x_2)\vee[(v_2-K)\vee 0]   	
 \]
and   for any $(y^0,y,\nu)\in W^{1,1}([0, \bar S];\R^{1+n+1})$,  introduce  the payoff
\begin{equation} \label{funz_obiettivo}
{\mathcal J}(y^0,y,\nu):= \Phi \left(y^0(0),y(0),y^0(\bar S),y(\bar S),\nu(\bar S)\right) \vee \max_{s \in [0,\bar S]} h(y^0(s),y(s))   .
\end{equation}
Fixed  a sequence $(\varepsilon_i)_i$ such that $\varepsilon_i  \downarrow 0$,  for each $i \in \N$  we consider the fixed end-time   optimal control problem:
\begin{equation*}
\left(\hat P_i\right)
\begin{cases}
\qquad\qquad\qquad\qquad\qquad\text{minimize } \,\,\,{\mathcal J}(y^0,y,\nu) \\
\text{over }({y^0}, y,\nu) \in W^{1,1}([0,\bar S]), \,\,\, (\zeta,\omega,\alpha)\in L^1([0,\bar S]) \text{ satisfying} \\
\ds\frac{dy^0}{ds}(s) = (1+\zeta(s))(1-|\omega(s)|) \qquad \text{a.e. } s \in [0, \bar S]  \\ 
\ds\frac{dy}{ds}(s) = (1 + \zeta)\Big[ f(y^0 ,y , \alpha)(1-|\omega|) + \sum_{j=1}^m g_j(y^0,y) \omega^j \Big](s) \,\, \text{ a.e. } s \in [0,\bar S]  \\ 
\ds\frac{d\nu}{ds}(s) = (1+\zeta(s)) |\omega(s)| \qquad \text{ a.e. } s \in [0,\bar S]  \\ 
 \nu(0)=0   \\
\w(s)\in (1-\varepsilon_i)(\C \cap \B), \,\,\, \zeta(s) \in [-1/2,1/2], \,\,\, \alpha(s) \in A \quad \text{ a.e. } s \in [0,\bar S].
\end{cases}
\end{equation*}
We will call an element $(\zeta,\omega,\alpha, {y^0}, y,\nu)$ satisfying the constraints in $(\hat P_i)$ a {\em feasible process for problem $(\hat P_i)$.} 
For every $i \in \N$, let   $(\bar S,  \hat\w_i^0,\hat{\w_i},\hat\alpha_i, \hat y_i^0, \hat{y_i},\hat{\nu_i})$ be the extended sense process  in which  $(\hat y^0_{i},\hat y_{i},  \hat \nu_{i} )(0)  = (\bar y^0, \bar y, \bar\nu)(0)$
and 
\[ 
\left(\hat\w_i^0 , \hat{\w_i},\hat\alpha_i \right)  (s) \,:=\left\{
 \begin{array}{l} \left(\varepsilon_i, (1- \varepsilon_i) \frac{\bar\w(s)}{|\bar\w(s)|},\bar\alpha(s) \right) \qquad \text{if $ \bar\w^0(s) <\varepsilon_i$} \\ \, \\
 \left( \bar\w^0(s) ,  {\bar\w}(s),\bar\alpha(s) \right)  \qquad\qquad \text{ if $\ {\bar\w}^0(s) \geq \varepsilon_i$. }
 \end{array} \right.
 \]
Notice that  $(\bar S,  \hat\w_i^0,\hat{\w_i},\hat\alpha_i, \hat y_i^0, \hat{y_i},\hat{\nu_i})$ is an embedded strict sense process for the extended problem \eqref{Pe}, since ${\hat\w}_i^0(s) \geq \varepsilon_i >0$ for a.e.  $s \in [0, \bar S]$. Moreover,  $(0, {\hat\w}_i,\hat\alpha_i , \hat y_i^0,\hat y_i, \hat\nu_i)$ is a feasible  process for $(\hat P_i)$, since ${\hat\w}_i(s) \in (1-\varepsilon_i)(\C \cap \B)$ a.e. $s\in[0,\bar S]$.  
Furthermore, 
\begin{equation} \label{conv_phi}
\hat\alpha_i\equiv \bar\alpha, \qquad  \left\|\left( {\hat\w}_i^0, {\hat\w}_i \right)- \left(  {\bar\w}^0,  {\bar\w}^0 \right) \right\|_{L^{\infty}([0,\bar S])} \to 0 \qquad \text{as $i \to \infty$},
\end{equation}
therefore, by the continuity of the input-output map $(\zeta,\omega,\alpha)\mapsto (y^0,y,\nu)$,  we have:
\begin{equation} \label{conv_haty}
\left\|\left(\hat y_i^0,\hat y_i , \hat\nu_i\right)- \left( \bar y^0, \bar y, \bar\nu \right) \right\|_{L^{\infty}([0,\bar S])} \to 0 \qquad \text{as $i \to \infty$}.
\end{equation}
Since ${\mathcal J}$ is nonnegative and  vanishes at $(\bar S, \bar\w^0 ,\bar\w,\bar\alpha , \bar y^0 ,\bar y , \bar\nu )$, by  the $L$-Lipschitz continuity of $h$ and the $1$-Lipschitz continuity of $d_{\T_0}(\cdot)$,  (\ref{conv_haty})    implies that  there exist a sequence $\rho_i \downarrow 0$ such that, for every $i\in\N$,  $(0, {\hat\w}_i,\hat\alpha_i , \hat y_i^0,\hat y_i, \hat\nu_i)$ has cost not greater than $\rho_i^2$, namely  is a $\rho_i^2$-minimizer for the problem $( \hat P_i )$.

\vsm
{\it Step 2.} If we endow the set of feasible processes for problem $(\hat P_i)$, say $\hat\Gamma_i$,  with the  
distance
\[
\begin{split}
& \hat d\left((\zeta,\omega,\alpha, {y^0}, y,\nu)\,,\, (\tilde\zeta,\tilde\omega,\tilde\alpha, \tilde y^0, \tilde y,\tilde\nu) \right) :=|(y^0,y)(0) - (\tilde y^0,\tilde y)(0)| \\
&\qquad\qquad   + \|\omega - \tilde\omega\|_{L^1([0,\bar S])} + \ell \{ s\in[0,\bar S]\,:\,(\zeta(s), \alpha(s))\neq (\tilde\zeta(s),\tilde\alpha(s))   \}   
\end{split}
\]
for every pair $(\zeta,\omega,\alpha, {y^0}, y,\nu)$, $(\tilde\zeta,\tilde\omega,\tilde\alpha, \tilde y^0, \tilde y,\tilde\nu)\in \hat\Gamma_i$,  $( \hat P_i )$ can be seen as an optimization problem with continuous cost over the complete metric space   $(\hat\Gamma_i,\hat d)$.
Hence, by Ekeland's Principle, if we introduce the function 
$$
\ell_i(s,\zeta, a):= \chi_{\{(\zeta,a) \neq (\zeta_i(s), \alpha_i(s))\}} \qquad \forall (s,\zeta,a)\in[0,\bar S]\times\left[-\frac{1}{2},\frac{1}{2}\right]\times A,
$$
 for any $i \in \N$ there is a feasible process $(\zeta_i,\omega_i,\alpha_i, y_i^0,y_i,\nu_i)$  for $( \hat P_i )$ which is a minimizer of
 \begin{equation*}
\left(P_i\right)
\begin{cases}
\text{minimize} \,\, {\mathcal J}(y^0,y,\nu) +\rho_i \Bigl(  \left| (y^0,y)(0) - (y_i^0,y_i)(0)   \right| \\ 
\qquad\qquad\qquad\qquad\qquad+ \int_0^{\bar S}  \left[ |\omega - \omega_i|(s) + \ell_i(s,\zeta(s), \alpha(s)) \right] \,ds  \Bigr)   \\
 \text{over }({y^0}, y,\nu) \in W^{1,1}([0, \bar S]), \,\,\, (\zeta,\omega,\alpha)\in L^1([0,\bar S]) \text{ satisfying} \\
\ds\frac{dy^0}{ds}(s) = (1+\zeta(s))(1-|\omega(s)|) \qquad \text{a.e. } s \in [0, \bar S]  \\ 
\ds\frac{dy}{ds}(s) = (1 + \zeta)\Big[ f(y^0,y, \alpha)(1-|\omega|) + \sum_{j=1}^m g_j(y^0,y) \omega^j   \Big](s)  \,\, \text{ a.e. } s \in [0,\bar S]  \\ 
\ds\frac{d\nu}{ds}(s) = (1+\zeta(s)) |\omega(s)| \qquad \text{ a.e. } s \in [0,\bar S]  \\ 
 \nu(0)=0   \\
 \zeta(s) \in [-1/2,1/2], \,\,\, \w(s)\in (1-\varepsilon_i)(\C \cap \B), \,\,\, \alpha(s) \in A   \,\,\, \text{ a.e. } s \in [0,\bar S] 
\end{cases}
\end{equation*}
and verifies 
\begin{equation} \label{ek}
 \ds \left| \left( \hat y_i^0,\hat y_i  \right)(0)-\left( y_i^0,y_i  \right)(0)  \right| + \int_0^{\bar S} \left[ |\hat\omega_i - \omega_i|(s) + \ell_i(s, 0, \bar\alpha(s)) \right] \,ds \leq \rho_i \to 0. 
\end{equation}
Thus, by (\ref{conv_phi}) and (\ref{ek}) it follows that, as $i \to +\infty$,  
 \begin{equation} \label{conv_y}
 \left\| \left( y_i^0, y_i,\nu_i  \right)-\left( \bar y^0, \bar y, \bar\nu  \right)  \right\|_{L^{\infty}([0,\bar S])} \to 0, \qquad 
\left\| \w_i-  \bar \w  \right\|_{L^1([0,\bar S])} \to 0 
\end{equation}
\begin{equation} \label{conv_mis}
 \ell\{ s\in[0,\bar S]\,:\,(\zeta_i(s),\alpha_i(s)) \neq (0,\bar\alpha(s))  \} \to 0 
\end{equation}
 so that, eventually passing to a subsequence,  $(\w_i)_i$ converges to $\bar \w$ almost everywhere. 
 
   Let us now show that,  through  suitable reparameterization techniques, the sequence of minimizing processes  $(\zeta_i,\omega_i,\alpha_i, y_i^0,y_i,\nu_i)$ can be associated to a sequence of embedded strict processes converging to the original  isolated process $(\bar S, \bar\w^0 ,\bar\w,\bar\alpha , \bar y^0 ,\bar y , \bar\nu )$.

\noindent Precisely, for each $i\in \N$, let us consider the surjective,  bi-Lipschitz continuous, and strictly increasing function $\sigma_i:[0,\bar S]\to [0,\tilde S_i]$, given by 
 $$
 \ds\sigma_{i}(s):= \int_0^s (1+ \zeta_{i}(r))\text{ d}r , \qquad \tilde S_i:= \sigma_i(\bar S).
 $$
Using as reparameterization the  inverse function $\sigma^{-1}_i:[0,\tilde S_i]\to [0,\bar S]$,  we derive that the corresponding process $(\tilde S_i, \tilde\w^0_{i}, \tilde\w_{i},\tilde \alpha_i, \tilde y^0_{i} ,\tilde y_{i}, \tilde \nu_i )$, where
\begin{equation} \label{maggiore}
( \tilde y^0_{i} ,\tilde y_{i}, \tilde \nu_i ):=(y^0_{i}, y_{i}, \nu_i,)\circ\sigma_{i}^{-1}, \quad  (\tilde\w^0_{i}, \tilde\w_{i},\tilde \alpha_i):= (1-|\w_i|, \w_i, \alpha_i)\circ\sigma_{i}^{-1},
\end{equation}
is an embedded  strict sense process for problem \eqref{Pe}.  In particular,   we have
\begin{equation} \label{list2}
\begin{split}
&{\tilde\w}_i^0(s)=1- \left| {\tilde\w}_i (s)\right|, \qquad   {\tilde\w}_i (s) \in (1-\epsilon_{i})\,(\C \cap \B_m) \ \mbox{ a.e. } s \in [0,\tilde S_i], \\
&(\tilde y^0_i(0), \tilde y_{i}(0), \tilde y^0_{i}(\tilde S_{i}), \tilde y_{i}(\tilde S_i),  \tilde \nu_i(\tilde S_{i}))\;= \;(y^0_i(0), y_{i}(0), y^0_{i}(\bar S),  y_{i}(\bar S), \nu_i(\bar S))\,.
\end{split}
\end{equation}
Hence,  we deduce from  (\ref{conv_y}) that, for $i$ sufficiently large, 
\begin{equation}
\label{list3}
d_\infty\Big((\tilde y^0_i(0),\tilde y^0_i(\tilde S_i), \tilde y_i,\tilde\nu_i), (\bar y^0(0),\bar y^0(\bar S),\bar y , \bar \nu)\Big)<\delta,
\end{equation}
where $\delta>0$ is the constant appearing in  Definition \ref{isolated}, with reference to the isolated feasible  extended sense process $(\bar S, \bar\w^0 ,\bar\w,\bar\alpha , \bar y^0 ,\bar y , \bar\nu )$.
As a consequence, for all $i$ large enough, $(\tilde S_i, \tilde\w^0_{i}, \tilde\w_{i},\tilde \alpha_i, \tilde y^0_{i} ,\tilde y_{i}, \tilde \nu_i )$ cannot be a feasible embedded strict sense process, namely, it must violate either the endpoint constraints or the state constraint. By  (\ref{maggiore}), (\ref{list2}) this implies  that ${\mathcal J}( y_i^0,y_i,\nu_i)>0$, namely, at least one of the following three inequalities holds true:
\begin{equation} \label{violazione} 
d_{\T_0}(y_i^0(0),y_i(0),y_i^0(\bar S),y_i(\bar S))>0, \quad \nu_i(\bar S)>K, \quad \max_{s \in [0,\bar S]} \,  h(y_i^0(s),y_i(s))  >0.
\end{equation}
 In the following, as is clearly not restrictive, we will always assume that the properties valid from a certain index onwards, apply to each index   $i\in\N$.

\vsm
{\it Step 3.}  For each $i\in\N$, define $c_i := \max_{s \in [0,\bar S]} \, h(y_i^0(s),y_i(s))$ and set
$$
\tilde h(t,x,c) := h(t,x) - c \qquad \forall (t,x,c)\in\R^{1+n+1}.
$$
The process  $(\zeta_i,\omega_i, \alpha_i, y_i^0,y_i,\nu_i,c_i)$ turns out to be a minimizer for 
\begin{equation*}
\left(Q_i\right)
\begin{cases}
\text{minimize} \,\,\,\Phi(y^0(0),y(0),y^0(\bar S),y(\bar S),\nu(\bar S)) \vee c(\bar S)  \\
\,\,\,\,\,  + \rho_i \left(  \left| (y^0,y)(0) - (y_i^0,y_i)(0)   \right| +
 \int_0^{\bar S} \left[ |\omega - \omega_i|(s) + \ell_i(s, \zeta(s),\alpha(s)) \right]\,ds  \right) \\
 \text{over }({y^0}, y,\nu,c) \in W^{1,1}([0, \bar S]), \,\,\, (\zeta,\omega,\alpha)\in L^1([0,\bar S]) \text{ satisfying} \\
\ds\frac{dy^0}{ds}(s) = (1+\zeta(s))(1-|\omega(s)|) \qquad \text{a.e. } s \in [0, \bar S]  \\ 
\ds\frac{dy}{ds}(s) = (1 + \zeta)\Big[ f(y^0,y, \alpha)(1-|\omega|) + \sum_{j=1}^m g_j(y^0,y) \omega^j   \Big](s) \,\,\text{ a.e. } s \in [0,\bar S]  \\ 
\ds\frac{d\nu}{ds}(s) = (1+\zeta(s)) |\omega(s)| \qquad \text{ a.e. } s \in [0,\bar S]  \\[1.5ex] 
\ds\frac{dc}{ds}(s) =0  \qquad \text{ a.e. } s \in [0,\bar S]  \\
 \nu(0)=0   \\
 \tilde h(y^0(s),y(s),c(s))  \leq 0 \qquad \text{for all } s\in[0,\bar S]   \\
(\zeta,\omega,\alpha,c)(s)\in [-1/2,1/2]\times(1-\varepsilon_i)(\C \cap \B) \times A \times \R  \quad \text{ a.e. } s \in [0,\bar S].
\end{cases}
\end{equation*}
Our aim is now to  apply the Pontryagin Maximum Principle to problem $(Q_i)$ with reference to the minimizer $(\zeta_i,\omega_i, \alpha_i, y_i^0,y_i,\nu_i,c_i)$. Preliminarily, let us observe that, passing eventually to a subsequence,   we may  assume that either $c_i > 0$ for each $i \in \N$ or $c_i \leq0$ for each $i \in \N$.

\vsm
 
Assume first that $c_i > 0$ for each $i \in \N$. 
  Fix $i\in\N$ and set $\w_{0_i}:=1-|\w_i|$.   In the Maximum Principle, several generalized subdifferentials are involved which it is convenient to make as explicit as possible.  First of all, the condition `$h(y_i^0(s),y_i(s)) - c_i >0$' implies `$h(y_i^0(s),y_i(s))>c_i >0$', so that 
$
\partial_{t,x,c}^{^{>}} \tilde h(y_i^0(s),y_i(s),c_i) = \partial_{t,x}^{^{>}} h(y_i^0(s),y_i(s)) \times \{ -1\}.
$
Furthermore,   by the  `max rule'  of subdifferential calculus  (see \cite[Th. 5.5.2]{OptV}), the properties of the subdifferential of the distance function (see  \cite[Lemma 4.8.3]{OptV}), and  \eqref{violazione}, we have that $(\gamma_{0_i}^1, \gamma_{i}^1,\gamma_{0_i}^2,\gamma_{i}^2,\gamma_{\nu_i}^2,\gamma_{c_i}^2)\in \partial \left(\Phi(y^0(0),y(0),y^0(\bar S),y(\bar S),\nu(\bar S)) \vee c(\bar S) \right)$ implies that there are some $\sigma^1_i$, $\sigma_i^2$, $\sigma_i^3\ge0$ with $\sum_{k=1}^3\sigma_i^k=1$, such that
$$
\begin{array}{l}
 \ds(\gamma_{0_i}^1, \gamma_{i}^1,\gamma_{0_i}^2,\gamma_{i}^2)\in\sigma_i^1\, \left(\partial d_{\T_0}(y_i^0(0),y_i(0),y_i^0(\bar S),y_i(\bar S))\cap \partial\B_{1+n+1+n}\right), \\[1.5ex]
\quad\ds\gamma_{\nu_i}^2= \sigma_i^2 \quad \text{(since $\partial \left((\nu_i(\bar S)-K)\vee 0\right)=1$ when $\nu_i(\bar S)>K$),} \qquad 
  \gamma_{c_i}^2=\sigma_i^3,
 \end{array}
$$
and $\sigma_i^k=0$ when the maximum in $\Phi \left(y_i^0(0),y_i(0),y_i^0(\bar S),y_i(\bar S),\nu_i(\bar S)\right) \vee c_i(\bar S)$ is strictly greater than the $k$-th term in the maximization.
Thus,  the  Maximum Principle in \cite[Th. 9.3.1]{OptV} yields the existence of some multipliers $(p_{0_i},p_i,\pi_i,r_i) \in W^{1,1}([0, \bar S]; \R^{1+n+1+1})$  associated with  $(y_i^0, y_i,\nu_i, c_i)$,  $\mu_i \in NBV^{+}([0,\bar S]; \R)$, $\lambda_i\geq0$, $\sigma^1_i$, $\sigma_i^2$, $\sigma_i^3\ge0$ with $\sum_{k=1}^3\sigma_i^k=1$, and Borel-measurable, $\mu_i$-integrable functions $(m_{0_i},m_i): [0,\bar S] \to \R^{1+n}$, such that:
\begin{itemize}
\item[(i)$'$] $\| p_{0_i} \|_{L^{\infty}} +\| p_i \|_{L^{\infty}}+ \| \mu_i \|_{TV} + \lambda_i + \| r_i \|_{L^{\infty}} + \| \pi_i \|_{L^{\infty}} =1$;

\item[(ii)$'$]  $\ds\left(-\frac{dp_{0_i}}{ds}(s), -\frac{dp_i}{ds}(s)\right) \in co \text{ } \partial_{t,x} \Big\{ \Big[q_i(s) \cdot \big(f((y_i^0, y_i,\alpha_i)(s))\omega_i^0(s)$\\
\textcolor{white}{r} \quad\qquad\qquad\qquad\qquad\qquad\qquad $+\sum_{j=1}^{m}g_{j}((y_i^0, y_i)(s))\omega_i^j(s)  \big)\Big] (1+\zeta_i(s)) \Big\}$ \\ 
and $d\pi_i(s)/ds = dr_i(s)/ds =0$ for a.e.  $s \in [0, \bar S]$;

\item[(iii)$'$] $\left(p_{0_i}(0), p_i(0), -q_{0_i}(\bar S), -q_i(\bar S)  \right) \in  \lambda_i \rho_i \B_{1+n} \times \{0_{1+n} \}$ \\ 
\textcolor{white}{sdfkssdfsddjf}$+\lambda_i \sigma_i^1  \left(\partial   d_{\T_0}(y_i^0(0),y_i(0),y_i^0(\bar S),y_i(\bar S))\cap \partial\B_{1+n+1+n}\right)$, \\
 $-\pi_i=\lambda_i \,\sigma_i^2$, \ \
 $r_i(0)=0$, \ \  $-r_i(\bar S) + \int_{[0,\bar S]} \mu_i(ds)=\lambda_i \sigma_i^3$; 

\item[(iv)$'$] $\int_0^{\bar S} \{ [ q_i\cdot( f(y_i^0,y_i,\alpha_i) \omega_i^0 + \sum_{j=1}^{m}g_{j}(y_i^0,y_i)\omega_i^j ) + q_{0_i} \omega_i^0  + \pi_i |\omega_i| ] (1+\zeta_i)\}\text{d}s$ \\
\textcolor{white}{e}$\geq \int_0^{\bar S} \{ [ q_i \cdot ( f(y_i^0,y_i,\alpha) \omega^0 + \sum_{j=1}^{m}g_{j}(y_i^0,y_i)\omega^j ) + q_{0_i} \omega^0  + \pi_i |\omega| ] (1+\zeta) $ \\ 
\textcolor{white}{e}\qquad\qquad\qquad\qquad\qquad  $- \lambda_i \rho_i \left[ |\omega - \omega_i| + \ell_i(s, \zeta(s),\alpha(s)) \right] \}\text{ d}s$  \\
for all measurable selectors $(\zeta,\omega,\alpha)$ of $[-\frac{1}{2},\frac{1}{2}]\times(1-\varepsilon_i)(\C \cap \B_m)\times A$;

\item[(v)$'$] $(m_{0_i},m_i)(s) \in \partial_{t,x}^{^{>}} h\left(y_i^0(s),  y_i(s)\right)$ \qquad $\mu_i$-a.e. $s \in [0, \bar S]$,

\item[(vi)$'$] $spt(\mu_i) \subseteq \{ s \text{ : } h\left( y_i^0(s), y_i(s)\right) - c_i = 0 \}$,
\end{itemize}
where
\[
(q_{0_i},q_i)(s) := 
\begin{cases}
\ds (p_{0_i},p_i)(s)+ \int_{[0,s[} (m_{0_i},m_i)(\tau) \mu_i(d\tau) \qquad\,\, s\in[0, \bar S[, \\
\ds(p_{0_i}, p_i)(\bar S) + \int_{[0, \bar S]} (m_{0_i},m_i)(\tau) \mu_i(d\tau) \qquad s=\bar S.
\end{cases}
\]
Observe that, for each $i \in \N$,  by (ii)$'$ and (iii)$'$  we have  
\begin{equation}\label{rmu}
r_i\equiv0, \qquad  \| \mu_i \|_{TV}= \int_{[0,\bar S]} \mu_i(ds)=\lambda_i \sigma_i^3, \qquad |\pi_i|=\lambda_i \,\sigma_i^2;
\end{equation}
 furthermore,   $\| (m_{0_i},m_i) \|_{L^{\infty}} \leq L$ by (v)$'$ and the $L$-Lipschitz continuity of $h$. 
Then  by (iii)$'$ and \eqref{rmu}, we get
\[
\begin{split}
\lambda_i(1-\sigma_i^3) -\lambda_i \rho_i&\leq \left| \left( p_{0_i}(0), p_i(0), -q_{0_i}(\bar S), -q_i(\bar S)   \right)  \right|+| \pi_i| \\ 
&\leq 2 \| p_{0_i} \|_{L^{\infty}} +2\| p_i \|_{L^{\infty}}+ 2L \| \mu_i \|_{TV} + | \pi_i |.
\end{split}
\] 
By this  estimate, \eqref{rmu},  the non-triviality condition (i)$'$,     
 and using the facts that $\rho_i \leq \frac{1}{2}$ for $i$ sufficiently large and $\lambda_i \in  [0,1]$, for such $i$ we get
\[	
\begin{array}{l} 3 \| p_{0_i} \|_{L^{\infty}} +3\| p_i \|_{L^{\infty}}+ (2L+2) \| \mu_i \|_{TV} +3 | \pi_i |  \\[1.5ex]
\qquad\qquad \geq  \lambda_i(1-\sigma_i^3) -\lambda_i \rho_i+1-\lambda_i+\lambda_i\sigma_i^2+\lambda_i\sigma_i^3 =1+\lambda_i\sigma_i^2-\lambda_i \rho_i\ge 
\frac{1}{2}.
\end{array}	
\]
Hence,    scaling the multipliers, we obtain  
\begin{equation}\label{stima}
 \| p_{0_i} \|_{L^{\infty}} +\| p_i \|_{L^{\infty}}+ \| \mu_i \|_{TV} +  | \pi_i | = 1, \qquad \lambda_i \leq \tilde L:=6\vee 4(1+L).
\end{equation}


\vsm
 
Suppose now $c_i \leq 0$ for each $i \in \N$.  In this case, by \eqref{violazione}, either  $\nu_i(\bar S)>K$ or $d_{\T_0}(y_i^0(0),y_i(0),y_i^0(\bar S),y_i(\bar S))>0$. Thus, for $\varepsilon >0$ suitably small,  the process  $(\zeta_i,\omega_i, \alpha_i, y_i^0,y_i,\nu_i,\hat c_i)$ with $\hat c_i := c_i + \varepsilon$ is still a minimizer for problem $(Q_i)$ and, in addition, it verifies  $h(y_i^0(s),y_i(s)) -\hat c_i <0$  for all  $s\in[0,\bar S]$ (namely, the state constraint is inactive on $[0,\bar S]$). 
  Hence, by applying the Maximum Principle for problem $(Q_i)$ with reference to this minimizer we deduce the existence of multipliers  $(p_{0_i},p_i,\pi_i,r_i) \in W^{1,1}([0, \bar S]; \R^{1+n+1+1})$,  which satisfy conditions (i)$'$--(vi)$'$ with $\mu_i=0$, $\sigma_i^3=0$, $\pi_i\le0$, and $\lambda_i > 0$.
In this case,  by considering  again $i$ sufficiently large to have  $\rho_i \leq \frac{1}{2}$, from  (iii)$'$ we get
\[ \lambda_i(1-\rho_i) \leq \left| \left( p_{0_i}(0), p_i(0), -q_{0_i}(\bar S), -q_i(\bar S) \right)\right|+| \pi_i| \leq 	2 \| p_{0_i} \|_{L^{\infty}} +2\| p_i \|_{L^{\infty}}  + | \pi_i |,\]
and,  scaling the multipliers appropriately after summing (i)$'$, we finally obtain 
\begin{equation} \label{stima2}
 \| p_{0_i} \|_{L^{\infty}} +\| p_i \|_{L^{\infty}}+ | \pi_i | = 1, \qquad \lambda_i \leq 6 \quad ( \leq \tilde L).
\end{equation}

\vsm
{\it Step 4.} From the previous step,  we arrive at the following properties (for either the case where   $c_i > 0$ for each $i \in \N$ or the case where $c_i \leq0$ for each $i \in \N$): for any $i\in\N$, 
 there exist $(p_{0_i},p_i) \in W^{1,1}([0,\bar S]; \R^{1+n})$, $\pi_i\le0$,  $\mu_i \in NBV^{+}([0,\bar S]; \R)$ and Borel-measurable,   $\mu_i$-integrable functions $(m_{0_i},m_i): [0,\bar S] \to \R^{1+n}$,  such that:
\begin{itemize}
\item[ (i)] $\| p_{0_i} \|_{L^{\infty}} +\| p_i \|_{L^{\infty}}+ \| \mu_i \|_{TV} + | \pi_i |=1$,

\item[(ii)]  $\ds\left(-\frac{dp_{0_i}}{ds}(s), -\frac{dp_i}{ds}(s)\right) \in co \text{ } \partial_{t,x} \Big\{ \Big[q_i(s) \cdot \big(f(y_i^0(s), y_i(s),\alpha_i(s))\omega_i^0(s)$\\
\textcolor{white}{r}  \qquad\qquad $+\sum_{j=1}^{m}g_{j}(y_i^0(s), y_i(s))\omega_i^j(s)  \big)\Big] (1+\zeta_i(s)) \Big\}$  
  for a.e.  $s \in [0, \bar S]$;

\item[(iii)]  $\left(p_{0_i}(0), p_i(0), -q_{0_i}(\bar S), -q_i(\bar S), -\pi_i  \right)$ \\ 
\textcolor{white}{r}  \quad$\in[0,\tilde L]\cdot\partial   \Phi(y_i^0(0),y_i(0),y_i^0(\bar S),y_i(\bar S),\nu_i(\bar S)) +\tilde L \rho_i \B_{1+n} \times \{0_{1+n} \}\times\{0\}$; 

\item[(iv)] $\int_0^{\bar S} \{ [ q_i \cdot ( f(y_i^0,y_i,\alpha_i) \omega_i^0 + \sum_{j=1}^{m}g_{j}(y_i^0,y_i)\omega_i^j ) + q_{0_i} \omega_i^0  + \pi_i |\omega_i| ] (1+\zeta_i)\}\text{d}s$ \\
\textcolor{white}{r} \qquad\qquad  $\geq \int_0^{\bar S} \{ [ q_i \cdot ( f(y_i^0,y_i,\alpha) \omega^0 + \sum_{j=1}^{m}g_{j}(y_i^0,y_i)\omega^j )$\\
 \textcolor{white}{r} \qquad\qquad \qquad\qquad $+ q_{0_i} \omega^0  + \pi_i |\omega| ] (1+\zeta) - 3\tilde L \rho_i \}\text{d}s$, \\
for all measurable selectors $(\zeta,\omega,\alpha)$ of $[-\frac{1}{2},\frac{1}{2}]\times(1-\varepsilon_i)(\C \cap \B_m)\times A$;

\item[(v)] $(m_{0_i},m_i)(s) \in \partial_{t,x}^{>} h\left(y_i^0(s),  y_i(s)\right)$ \qquad $\mu_i$-a.e. $s \in [0, \bar S]$,

\item[(vi)] $spt(\mu_i) \subseteq \{ s \text{ : } h\left( y_i^0(s), y_i(s)\right) - c_i = 0 \}$,
\end{itemize}
where
\[
(q_{0_i},q_i)(s) := 
\begin{cases}
\ds(p_{0_i},p_i)(s) + \int_{[0,s[} (m_{0_i},m_i)(\tau) \mu_i(d\tau) \qquad\,\, s\in[0, \bar S[, \\
\ds(p_{0_i},p_i)(\bar S) + \int_{[0, \bar S]} (m_{0_i},m)(\tau) \mu_i(d\tau)\qquad\,\, s=\bar S.
\end{cases}
\]
Here $\tilde L$ is the same constant as in \eqref{stima}.  By Banach-Alaoglu's Theorem, there exist a subsequence of $(\mu_i)_i$,  $ \mu \in NBV^{+}([0,\bar S]; \R)$,  $(m_0,m):[0,\bar S] \to \R\times\R^n$ Borel measurable and $\mu$-integrable,  such that $\mu_i\overset{*}{\rightharpoonup} \mu$ weakly* in $C^*([0,\bar S])$ and  $m_i \mu_i(ds)\overset{*}{\rightharpoonup} m \mu(ds)$,  $m_{0_i} \mu_i(ds) \overset{*}{\rightharpoonup}  m_0 \mu(ds)$ (see  \cite[Proposition 9.2.1]{OptV}). Furthermore, 
$\int_{[0,s[} (m_{0_i},m_i)(\tau) \mu_i(d\tau) \to \int_{[0,s[} (m_0,m)(\tau) \mu(d\tau)$ for a.e.  $s \in [0,\bar S]$;  the  real sequence $(\pi_i)$ is bounded;    the  functions $(p_{0_i}, p_i)$ are uniformly bounded and  have uniformly integrable, bounded derivatives. Hence,  there exist $\pi\le0$ and $(p_0,p) \in W^{1,1}([0,\bar S]; \R^{1+n})$  (see e.g. \cite[Th 2.5.3 and Ch. 9]{OptV}) such that, eventually for a further subsequence,  $\pi_i \to \pi$, $(p_{0_i},p_i) \to (p_0,p)$   in $L^\infty$, and $\left(\frac{dp_{0_i}}{ds},\frac{dp_i}{ds}\right) \rightharpoonup \left(\frac{dp_0}{ds},\frac{dp}{ds}\right)$  weakly in $L^1$, as $i\to+\infty$.
By this   analysis it also follows  that the functions $(q_{0_i},q_i)$ are uniformly integrably bounded and verify for a.e. $s \in [0,\bar S]$,
\[
\ds \lim_{i\to+\infty}(q_{0_i},q_i)(s)= (q_0,q)(s):= (p_0,p)(s) + \int_{[0,s[} (m_0,m)(\tau) \mu(d\tau).
  \]
Hence, by the dominated convergence theorem, one has
\begin{equation} \label{compact2}
(q_{0_i},q_i) \to (q_0, q)  \qquad \text{in $L^1([0,\bar S])$.}
  \end{equation}
\noindent Passing to the limit as $i\to+\infty$  and using (\ref{conv_y}), by (i),(v),  and (vi)  we get 
\begin{equation} \label{quasi}
\| p_0 \|_{L^{\infty}} +\| p \|_{L^{\infty}}+ \| \mu \|_{TV} + | \pi | =1,
\end{equation}
\begin{equation} \label{fine2}
(m_0,m)(s) \in \partial_{t,x}^{>} h\left(\bar y^0(s),  \bar y(s)\right) \qquad \text{$\mu$-a.e. } s \in [0, \bar S],
\end{equation}
\begin{equation} \label{fine3}
spt(\mu) \subseteq \{ s \text{ : } h\left(\bar y^0(s), \bar y(s)\right) = 0 \}.
\end{equation}
Furthermore,  using that $\left(\bar{y}^0(0),\bar{y}(0),\bar{y}^0(\bar S), \bar y(\bar S), \bar\nu(\bar S)\right) \in \T$,  the properties of distance function, and the `max-rule' for subdifferentials,  by  (iii)  we have
\begin{equation}\label{fine4}
\begin{array}{l}
 (p_0(0), p(0),  -q_0(\bar S),  -q(\bar S) , -\pi )   \in [0,\tilde L]\cdot \partial \Phi\left(\bar{y}^0(0),\bar{y}(0),\bar{y}^0(\bar S), \bar y(\bar S),\bar\nu(\bar S)\right) \\
\qquad\qquad  \subseteq [0,\tilde L]\cdot\left[ \left(N_{\T_0}\left(\bar{y}^0(0),\bar{y}(0),\bar{y}^0(\bar S), \bar y(\bar S), \bar\nu(\bar S)\right)  \cap \B_{1+n+1+n}  \right)\times\{0\} \right. \\
\qquad\qquad\qquad\qquad\qquad\qquad\qquad  \left. + \{0_{1+n+1+n}\}\times \partial  \left((\nu_i(\bar S)-K)\vee 0\right)\right]\\
\qquad\qquad\qquad \subseteq N_{\T_0}\left(\bar{y}^0(0),\bar{y}(0),\bar{y}^0(\bar S), \bar y(\bar S)\right) \times N_{]-\infty,K]}(\bar\nu(\bar S)), 
\end{array}
\end{equation}
Incidentally, from this relation we immediately deduce that   $\pi=0$ if $\bar\nu(\bar S) < K$. 
Passing to the limit in (iv), with the help of a measurable selection theorem, using (\ref{conv_y}), (\ref{conv_mis}) and the dominated convergence Theorem, we deduce that, for a.e. $s \in [0, \bar S]$,
\[ 
\begin{array}{l}
H\Big(\bar y^0(s), \bar y(s), q_0(s) , q(s),\pi ,\bar \w^0(s),\bar \w(s),\bar\alpha(s)\Big) \\
\ds \qquad\qquad=  \max_{(\zeta,\w^0,\w,a)\in  [-\frac{1}{2},\frac{1}{2}]\times\CC\times A}  H\Big(\bar y^0(s), \bar y(s), q_0(s) , q(s),\pi, \w^0,\w,a\Big)(1+\zeta).
 \end{array}
\]
Since 0 is in the interior of $[-\frac{1}{2},\frac{1}{2}]$, this implies that, for a.e.  $s \in [0, \bar S]$,
\begin{equation} \label{fine6}
\begin{array}{l}
H\Big(\bar y^0(s), \bar y(s), q_0(s) , q(s),\pi,  \bar \w^0(s),\bar \w(s),\bar\alpha(s)\Big) \\
\ds \quad\qquad\qquad\qquad\qquad=  \max_{(\w^0,\w,a)\in  \CC\times A}  H\Big(\bar y^0(s), \bar y(s), q_0(s) , q(s),\pi, \w^0,\w,a\Big)=0.
 \end{array}\end{equation}
To prove that $(p_0,p)$ verifies the adjoint equation in Theorem \ref{Th_pmp}, for each $i\in\N$  we set
\[	
\begin{array}{l}
\A_i := \Big\{ s\in [0,\bar S] \text{ : } \left(\zeta_i(s), \alpha_i(s)\right) =  \left(0,\bar\alpha(s)   \right)      \Big\} \subseteq [0,\bar S].
 \end{array}		\]
By  (\ref{conv_mis}) it follows that $\ell(\A_i) \to \bar S$  as  $i \to +\infty.$
Hence, by the $L$-boundedness of $f$ and $g_1,\dots,g_m$, we deduce that the functions $(y_i^0, y_i)$ are uniformly integrable, so that by \cite[Th. 2.5.3]{OptV} and (\ref{conv_y}), $\left(\frac{dy^0_i}{ds},\frac{dy_i}{ds}\right) \rightharpoonup \left(\frac{dy^0}{ds},\frac{dy}{ds}\right)$  weakly in $L^1$, as $i\to+\infty$.
Moreover,  for a.e. $s \in \A_i$, we have
 \[
\begin{array}{l}
\ds \left(\frac{dy^0_i}{ds},\frac{dy_i}{ds}\right)(s)= \Big( \w_i^0, f(y_i^0,y_i,\bar\alpha) \w^0_i+ \sum_{j=1}^m g_j(y_i^0,y_i)  \w^j_i \Big)(s),  \\
 \ds \left( -\frac{dp_{0_i}}{ds},-\frac{dp_i}{ds}\right)(s) \in co\, \partial_{t,x} \Bigl\{ q_i(s) \cdot \Big(f(y_i^0, y_i,\bar\alpha)  \w^0_i+ \sum_{j=1}^{m}g_{j}(y_i^0,y_i)  \w^j_i \Big)(s)\Bigr\},
\end{array}
\]
By (\ref{conv_y}), (\ref{compact2}), and   \cite[Theorem 2.5.3]{OptV}  we can conclude that, for a.e. $s \in [0, \bar S]$,\footnote{Notice that, since  the convex hull of a compact subset of some space $\R^k$ is compact, our hypotheses on the vector fields $f,g_1, \dots, g_m$  guarantee that the multifunction in the differential inclusion is closed. }  
\begin{equation}\label{fine7}
\left( -\frac{dp_{0}}{ds},-\frac{dp}{ds}    \right)(s)  \in co \, \partial_{t,x} \Bigl\{ q(s) \cdot \Big(f(\bar{y}^0,\bar y,\bar\alpha)\bar\w^0+ \sum_{j=1}^{m}g_{j}(\bar{y}^0,\bar y)\bar\w^j \Big)(s)\Bigr\}.
\end{equation}
In view of relations (\ref{fine2}), (\ref{fine3}), (\ref{fine4}),   (\ref{fine6}), and  (\ref{fine7}), to conclude the proof that the isolated feasible process $(\bar S, \bar\w^0 ,\bar\w,\bar\alpha , \bar y^0 ,\bar y , \bar\nu )$  is an abnormal extremal,    it remains only to show that  
\begin{equation} \label{fine1}
\| p_0 \|_{L^{\infty}} +\| p \|_{L^{\infty}}+ \| \mu \|_{TV}  \neq 0.
\end{equation}
Suppose by contradiction that   \eqref{fine1} is not true.  Then  $q_0 \equiv 0$, $q \equiv 0$ a.e.,   and by (\ref{quasi}) we deduce that  $\pi \neq 0$, which in turn implies $\bar\nu(\bar S)=K >0$. Thanks to these information and integrating (\ref{fine6}) in $[0,\bar S]$ we find that 
$0 = \int_0^{\bar S} \pi \, \left| \bar \w\right| \text{d}s = \pi \, \bar\nu(\bar S) = \pi \, K$, which is not possible.	 		 
\qed
\vsm
 
\begin{remark}[Multiple state constraints] In order to allow multiple state constraints $h_i$ for $i=1,\dots,N$, it is sufficient replacing the payoff  $\mathcal{J}$ in (\ref{funz_obiettivo}) by the function $
\Phi \left(y^0(0),y(0),y^0(\bar S),y(\bar S),\nu(\bar S)\right) \vee \max_{s \in [0,\bar S]} h_1(y^0(s),y(s)) \vee \dots \vee \max_{s \in [0,\bar S]} h_N(y^0(s),y(s)),
$
and making obvious changes to the preceding proof. \end{remark}

\section{Nondegeneracy, normality and no infimum gap}\label{S3}
As a  consequence of Theorem \ref{th_gap}, `normality implies no infimum gap'. Precisely,  as a corollary of the results in Section \ref{S2}, we have:
\begin{theorem}\label{Th_Norm}
 Assume hypotheses
 {\bf (H0)}-{\bf (H2)}   are satisfied.
 \begin{itemize}
 \item[{\rm (i)}] Suppose that there exists an  extended sense  minimizer for \eqref{Pe} which is a normal extremal.  Then there is no infimum gap.
\item[{\rm (ii)}]
Let $(\bar S, \bar\w^0 ,\bar\w,\bar\alpha , \bar y^0 ,\bar y , \bar\nu )$  be an extended sense  local minimizer  for the extended problem \eqref{Pe} which is a normal extremal. Then there is no local infimum gap at $(\bar S, \bar\w^0 ,\bar\w,\bar\alpha , \bar y^0 ,\bar y , \bar\nu )$.
\end{itemize}
  \end{theorem}   
 
As    observed in the Introduction, the above `normality test' is of more theoretical than practical interest (specially in the presence of state constraints). In this section we  identify some verifiable conditions guaranteeing that every set of multipliers is normal. To begin with, let us  introduce the notion  of {\em nondegenerate estremal}.

\begin{definition}[Nondegenerate Maximum Principle] \label{nondeg_def}
Given an extended sense local minimizer  for problem \eqref{Pe},  $(\bar S, \bar\w^0 ,\bar\w,\bar\alpha , \bar y^0 ,\bar y , \bar\nu )$, we say that the Maximum Principle is \textit{nondegenerate} when there is a choice of the multipliers $(p_0,p,\pi, \lambda, \mu)$ and of the functions $m_0$, $m$ that meets the conditions (i)--(vi) of Theorem \ref{Th_pmp}, and such that
\begin{equation} \label{nondeg}
\begin{aligned}
&\mu(]0,\bar S]) + \Vert q_0 \Vert_{L^{\infty}} +\Vert q \Vert_{L^{\infty}} + \lambda \neq 0 \qquad &\text{if } \qquad \bar y^0(\bar S) = \bar y^0(0), \\
&\mu(]0,\bar S]) + \Vert q \Vert_{L^{\infty}} + \lambda \neq 0 \qquad &\text{if } \qquad \bar y^0(\bar S) > \bar y^0(0),
\end{aligned}
\end{equation}
where $q_0$, $q$ are defined as in Theorem \ref{Th_pmp}.
\end{definition}

\begin{definition}[Nondegenerate normal and abnormal  extremals]
We say that a feasible extended sense process $(\bar S, \bar\w^0 ,\bar\w,\bar\alpha , \bar y^0 ,\bar y , \bar\nu )$  is an  {\em extremal of the nondegenerate Maximum Principle,}  in short,  {\it a nondegenerate extremal,}  if there exists a set of multipliers $(p_0,p,\pi, \lambda, \mu)$ and functions $m_0$,   $m$ which meet the conditions  of Theorem \ref{Th_pmp} and also satisfy   \eqref{nondeg}.  We call $(\bar S, \bar\w^0 ,\bar\w,\bar\alpha , \bar y^0 ,\bar y , \bar\nu )$   
a {\it nondegenerate normal extremal}   if all possible choices of multipliers as above  have $\lambda >0$, and  a \textit{nondegenerate abnormal extremal}  when  there exists at least one set of such multipliers with $\lambda =0$. 
 \end{definition}
As it is   easy to see,  a nondegenerate abnormal extremal is always an abnormal extremal, and, on the contrary, any normal extremal is also nondegenerate normal. 
%
 To obtain the converse implications,  we introduce  condition  \textbf{(CNa)} below.  In the following, we will often use the notation  
 $$
 \Omega := \{ (t,x) \,:\, h(t,x) \leq 0 \}.
 $$

\vsm
\noindent\textbf{Condition for nondegenerate abnormality}   \textbf{(CNa)}. {\em A   feasible extended sense process $(\bar S, \bar\w^0 ,\bar\w,\bar\alpha , \bar y^0 ,\bar y , \bar\nu )$ is said to verify condition  {\bf (CNa)}
 if 
 \begin{equation} \label{nonsmooth_condgeom}   
 \partial^{>}h (\bar y^0(0), \bar y(0))  \, \cap \,\left( -\text{{\rm proj}}_{(t_1,x_1)}(N_{\T_0}(\bar y^0(0), \bar y(0), \bar y^0(\bar S), \bar y(\bar S))) \right) =\emptyset.       	
\end{equation}}
%
%
\begin{remark}\label{RemNa} To clarify the geometrical meaning of condition   \textbf{(CNa)},  let us notice that, 
if  $(\bar y^0(0), \bar y(0))\in\text{Int}(\Omega)$, condition \eqref{nonsmooth_condgeom}  is trivially satisfied, since the hybrid subdifferential $\partial^{>}h (\bar y^0(0), \bar y(0))=\emptyset$.  Incidentally, observe that $\{(t,x) \,:\, h(t,x)<0\}\subseteq \text{Int}(\Omega)$ but the inclusion is in general strict. When instead  $(\bar y^0(0), \bar y(0))\in\partial\Omega$, \eqref{nonsmooth_condgeom} implies that  $0\notin \partial^{>}h (\bar y^0(0), \bar y(0))$. 
 If $h  \in C^2$ in a neighborhood of $(\bar y^0(0), \bar y(0))\in\partial\Omega$,  \eqref{nonsmooth_condgeom} simply  reads ($\nabla h (\bar y^0(0), \bar y(0))\ne0$ and)
\begin{equation}\label{smooth_condgeom}
\nabla h (\bar y^0(0), \bar y(0))\notin  -\text{{\rm proj}}_{(t_1,x_1)}(N_{\T_0}(\bar y^0(0), \bar y(0), \bar y^0(\bar S), \bar y(\bar S))). 
\end{equation}
Condition \eqref{smooth_condgeom}  is satisfied at $(\bar y^0(0), \bar y(0))$ with $h(\bar y^0(0), \bar y(0))=0$ and \linebreak $\nabla h (\bar y^0(0), \bar y(0))\ne0$,  when, for instance,   $\T_0=\T_1\times\T_2$  with $\T_1$, $\T_2$ closed subsets of $\R^{1+n}$,  $\T_1\subseteq\Omega$,  and  $N_{\T_1}(\bar y^0(0), \bar y(0))$ is pointed.\footnote{A cone ${\mathcal K}\subseteq\R^k$ is  {\it pointed} if it contains no line, i.e. if $\xi$,  $-\xi\in {\mathcal K}$ implies that $\xi=0$.}
 In this case, indeed,  \eqref{smooth_condgeom}  can be derived by the following relations
\[
\begin{split}
\partial^> h(\bar y^0(0), \bar y(0))&=\{\nabla h(\bar y^0(0), \bar y(0))\}\subseteq N^P_{\Omega}(\bar y^0(0), \bar y(0)) \\
&\subseteq N^P_{\T_1}(\bar y^0(0), \bar y(0))\subseteq N_{\T_1}(\bar y^0(0), \bar y(0)).
\end{split}
\] 
 \end{remark}
 
 \begin{remark}\label{RemNatindip} Consider the quite customary situation where  initial and final time are fixed and the state constraint is time independent, namely $\T_0=\{\bar t_1\}\times\T^1_0\times \T^2_0$ with $\T^1_0\subseteq\R^n$, $\T^2_0\subseteq\R^{1+n}$ closed subsets, and $h(t,x)=\bar h(x)$. In this case,   $\partial^>h(t,x)=\{0\}\times \partial^>\bar h(x)$ and  $N_{\{\bar t_1\}\times\T^1_0}(t,x)=\R\times N_{\T^1_0}(x)$ for all $(t,x)\in\R^{1+n}$. Hence, condition \eqref{nonsmooth_condgeom} reduces to 
$$
\partial^>\bar h(\bar y(0))\cap \left(-N_{\T^1_0}(\bar y(0))\right)=\emptyset.
$$
\end{remark}

\begin{proposition} \label{interno} 
Assume  {\bf (H0)}-{\bf (H2)} and suppose  that  $(\bar S, \bar\w^0 ,\bar\w,\bar\alpha , \bar y^0 ,\bar y , \bar\nu )$ is a   feasible extended sense process which is an abnormal  extremal, namely, there exist a set of multipliers $(p_0,p,\pi,\lambda,\mu)$,  and some functions   $(m_0,m)$  as in Theorem  \ref{Th_pmp}, with  $\lambda=0$. 
If    condition  {\rm\textbf{(CNa)}}   is satisfied, then  $(p_0,p,\pi,0,\mu)$,  $(m_0,m)$  verify  the strengthened non-triviality condition \eqref{nondeg}.
\end{proposition}

%

\noindent {\it Proof: }\, 
Assume that {\bf (CNa)}  is verified and  suppose by contradiction that
$
\Vert q_0 \Vert_{L^{\infty}} + \Vert q \Vert_{L^{\infty}} + \mu(]0, \bar S]) = 0.
$
Then, in view of  Theorem \ref{Th_pmp}, we have
\[	(p_0,p) \equiv - \mu(\{0\})(\xi_0,\xi),		\quad \mu(\{0\}) \neq 0, \quad (\xi_0,\xi) \in \partial^> h (\bar y^0(0), \bar y(0)).
\]
In particular,    $\mu(\{0\}) \neq 0$  implies that $(\bar y^0(0),\bar y(0)) \in\partial\Omega$.  By the transversality condition (iii) of  Theorem  \ref{Th_pmp}, it follows  that
\begin{equation} \label{intermedio2}
(\xi_0,\xi) \in \partial^> h (\bar y^0(0), \bar y(0))\cap\left( -\text{{\rm proj}}_{(t_1,x_1)}(N_{\T_0}(\bar y^0(0), \bar y(0), \bar y^0(\bar S), \bar y(\bar S)))\right),
\end{equation}
in contradiction with  \eqref{nonsmooth_condgeom}.  

To conclude it remains to show that $\Vert q \Vert_{L^{\infty}} + \mu(]0, \bar S])\ne0$ whenever $\bar y_0(\bar S)>\bar y_0(0)$. If we suppose by contradiction $\mu(]0,\bar S])+\|q\|_\infty = 0$, then Theorem \ref{Th_pmp}, (ii), (iv) and the first part of the proof     yield that $0\neq q_0$ is a constant and
 \bel{contradd}
 q_0\bar w_0(s)+\pi |\bar w(s)|=\max_{w_0\in[0,1]}\{\pi+(q_0-\pi)w_0\}=0 \quad \text{a.e. $s\in]0,\bar S[$.}
\eeq
 Then $\pi<0$ leads easily to  a contradiction. If $\pi=0$,  since $q_0\ne0$,   \eqref{contradd} yields $\bar w_0=0$ a.e., in contradiction with  $\int_0^{\bar S}\bar w_0(s)\,ds=\bar y^0(\bar S)-\bar y^0(0)>0$. 
\qed
\vsm

 As a straightforward consequence of Proposition  \ref{interno}, we have:
 \begin{proposition}\label{nondeg_norm}  Assume  {\bf (H0)}-{\bf (H2)}. Let $(\bar S, \bar\w^0 ,\bar\w,\bar\alpha , \bar y^0 ,\bar y , \bar\nu )$  be a feasible extended sense process verifying  condition {\bf (CNa)}. Then  $(\bar S, \bar\w^0 ,\bar\w,\bar\alpha , \bar y^0 ,\bar y , \bar\nu )$ is a nondegenerate normal extremal if and only if it is a  normal extremal.  
\end{proposition}  
 
\begin{remark}[Multiple state constraints]\label{RMSC2} The result of Proposition \ref{interno} can be easily extended to the case of  multiple state constraints $h_i$ for  $i=1,\dots, N$, by considering the  nondegeneracy condition  \eqref{nondeg}  with $(q_0,q)$ and $\mu$ as in Remark \ref{RMSC1} and by replacing the nondegenerate abnormality condition  \eqref{nonsmooth_condgeom} with
\begin{equation}\label{nondeg_multiple_constraints}
\begin{split}
&\left( \sum_{i=1}^N[0,+\infty[ \cdot\partial^{>}h_i (\bar y^0(0), \bar y(0))\right)  \,  \\
&\textcolor{white}{lskjdflsdkj} \cap \left( -\text{{\rm proj}}_{(t_1,x_1)}(N_{\T_0}(\bar y^0(0), \bar y(0), \bar y^0(\bar S), \bar y(\bar S))) \right) =\emptyset. \footnote{To simplify the notation, when  $X\subseteq\R^k$,  we use the convention  $X+\emptyset=X$. Thus, in particular,  $\emptyset+\emptyset=\emptyset$.}
\end{split}
\end{equation}
\end{remark}

We now  provide some sufficient conditions in the form of constraint and endpoint qualifications to  guarantee normality.  
In the following,  we will use the  notation 
$$
\ds F(t,x,\w_0,\w,a):=f(t,x,a)\w^0+\sum_{i=1}^mg_i(t,x)\w^i  \quad    \forall(t,x,\w^0,\w,a)\in \R^{1+n}\times\CC\times A,
$$
and, given a feasible extended sense process $(\bar S, \bar\w^0 ,\bar\w,\bar\alpha , \bar y^0 ,\bar y , \bar\nu )$ for problem \eqref{Pe}, we set  
$$
\ds \bar F(s):=F((\bar y^0 ,\bar  y, \bar \w^0,\bar \w,\bar\alpha)(s)) \quad \forall s\in[0,\bar S], \qquad  \bar z:= (\bar S, \bar\w^0 ,\bar\w,\bar\alpha , \bar y^0 ,\bar y , \bar\nu ). 
$$
  
 \noindent\textbf{Constraint qualifications for normality   \textbf{(CQn)$_b$},  \textbf{(CQn)$_f$}.}  {\em Let  $\bar z$  be a feasible extended sense process for the extended optimization problem \eqref{Pe}. 
 \begin{itemize}
\item[]\textbf{(CQn)$_b$}  We say that $\bar z$ meets condition {\bf(CQn)$_b$} if for every $s \in ]0,\bar S]$ where $(\bar y^0(s), \bar y(s))\in\partial\Omega$ there exist $\varepsilon$, $\delta > 0$ such that  
\[          \inf_{a \in A} \, \max_{(\xi_0,\xi) \in \partial^* h(\bar y^0(s), \bar y(s))}  \left[  \xi \cdot \left( f(\bar y^0(s), \bar y(s),a) - f(\bar y^0(s), \bar y(s), \bar\alpha(\sigma)) \right)  \right] < -\delta,    \]
\[       \inf_{\omega \in \C \cap \partial\B} \, \max_{(\xi_0,\xi) \in \partial^* h(\bar y^0(s), \bar y(s))}  \Big[  \xi \cdot \big( \sum_{j=1}^m g_j(\bar y^0(s), \bar y(s)) ( \omega^j - \frac{ {\bar\w}^j(\sigma)}{| {\bar\w}^j(\sigma)|} ) \big)  \Big] < -\delta  \,  \footnote{We use the convention that $\frac{\omega}{|\omega|} = 0$ when $\omega =0$.}       \]
for a.e. $\sigma \in  E(s,\varepsilon)$, defined as follows
\[E(s,\varepsilon):= \{  r \in [s-\varepsilon, s] \cap [0,\bar S]:\max_{(\xi_0,\xi) \in \partial^* h(\bar y^0(r), \bar y(r))} (\xi_0  {\bar\w}^0 + \xi \cdot \bar F)(r) \geq 0 \}  ;  \]
\item[]\textbf{(CQn)$_f$} We say that $\bar z$ meets condition {\bf(CQn)$_f$} if for every  $s \in [0,\bar S[$  where $(\bar y^0(s), \bar y(s))\in\partial\Omega$ there exist $\varepsilon$, $\delta > 0$ such that  
\[          \sup_{a \in A} \, \min_{(\xi_0,\xi) \in \partial^* h(\bar y^0(s), \bar y(s))}  \left[  \xi \cdot \left( f(\bar y^0(s), \bar y(s),a) - f(\bar y^0(s), \bar y(s), \bar\alpha(\sigma)) \right)  \right] > \delta,    \]
\[        \sup_{\omega \in \C \cap \partial\B} \min_{(\xi_0,\xi) \in \partial^* h(\bar y^0(s), \bar y(s))}  \left[  \xi \cdot ( \sum_{j=1}^m g_j(\bar y^0(s), \bar y(s)) ( \omega^j - \frac{ {\bar\w}^j(\sigma)}{|{\bar\w}^j(\sigma)|} ) )  \right] > \delta         \]
for a.e.  $\sigma \in  \Gamma(s,\varepsilon)$, defined as follows
\[
\Gamma(s,\varepsilon):= \{  r \in [s, s+\varepsilon] \cap [0,\bar S]:\min_{(\xi_0,\xi) \in \partial^* h(\bar y^0(r), \bar y(r))}  (\xi_0  {\bar\w}^0 + \xi \cdot \bar F)(r)  \leq 0 \} .
 \]
\end{itemize}}

\begin{remark}\label{Rem2cqn}  The  `forward' constraint qualification     {\bf (CQn)$_f$} is at our knowledge  new, while a version of the   `backward' constraint qualification    {\bf (CQn)$_b$} was first introduced in \cite{MS20}, as an adaptation to impulsive optimal control  of a   condition due to \cite{FoFr15}.
 In particular,  {\bf (CQn)$_b$},  {\bf (CQn)$_f$}  prescribe    that   drift and  fast dynamics satisfy  {\it separately}  outward or inward   pointing conditions, respectively. This is a stronger requirement than the existence of an inward/outward pointing velocity.  In fact, as discussed in \cite[Remark 4.4]{MS20}, the statement of Theorem \ref{Th_gapsuff} below  holds true even if  {\bf (CQn)$_b$},  {\bf (CQn)$_f$}  are  replaced, respectively,  with the (weaker) conditions: 
{\em \begin{itemize}
\item[]{\bf (CQn)$'_b$}   for every $s\in]0,\bar S]$  such that $(\bar y^0(s), \bar y(s))\in\partial\Omega$ there exist $\varepsilon$, $\delta>0$   and  a measurable control $(\tilde \w,\hat\alpha)$ taking values in $(\C\cap\partial \B)\times A$, satisfying    for all  $\sigma\in] s-\varepsilon, s[\cap[0,\bar S]$:
\bel{cqnbtot'}
\ds \max_{(\xi_0,\xi)\in\partial^*h(\bar y^0(\sigma),\bar y(\sigma))} \xi\cdot \Big(F(\bar y^0,\bar y,\bar \w^0,\hat \w,\hat\alpha)(s') -\bar F(s')\Big)<-\delta,  
\eeq
where $\hat \w:=(1-\bar \w^0)\tilde \w$, for a.e. $s'\in E(s,\varepsilon)$,  defined as in  {\bf (CQn)$_b$}; 

\item[]\textbf{(CQn)$'_f$} for every $s\in[0,\bar S[$  such that $(\bar y^0(s), \bar y(s))\in\partial\Omega$ there exist $\varepsilon$, $\delta>0$   and  a measurable control $(\tilde \w,\hat\alpha)$ taking values in $(\C\cap\partial \B)\times A$, satisfying    for all  $\sigma\in] s,s+\varepsilon[\cap[0,\bar S]$:
\begin{equation} \label{cqnftot'}
\min_{(\xi^0,\xi)\in\partial^*h(\bar y^0(\sigma),\bar y(\sigma))} \xi\cdot \Big(F(\bar y^0,\bar y,\bar \w^0,\hat \w,\hat\alpha)(s') -\bar F(s')\Big)>\delta,  
\end{equation}
where $\hat \w:=(1-\bar \w^0)\tilde w$, for a.e. $s'\in \Gamma(s,\varepsilon)$,  defined as in  {\bf (CQn)$_f$}.
\end{itemize}}
Notice that conditions \textbf{(CQn)$'_b$}, \textbf{(CQn)$'_f$} also cover  situations in which the drift does not depend on the ordinary control $a$, i.e. $f=f(t,x)$, unlike \textbf{(CQn)$_b$}, \textbf{(CQn)$_f$}, where in this case the hypotheses involving $f$  are clearly never fulfilled.
 \end{remark}

\begin{remark}\label{Ripfc} 
 The  constraint qualifications   {\bf (CQn)$_b$},  {\bf (CQn)$_f$}  are respectively inward/outward pointing  conditions at the boundary  which involve the minimizer but have to be satisfied just on a subset of instants at which the optimal trajectory has an outward/inward pointing velocity. As discussed in detail in  \cite[Lemma 4.1]{MS20},  when the constraint function $h$ is  smooth, {\bf (CQn)$_b$},  {\bf (CQn)$_f$}  could be replaced  by  the following simpler  conditions: 
{\em \begin{itemize}
\item[]{\bf (IPFCn)}$_b$ for every $s\in ]0,\bar S]$ such that $(\bar y^0(s), \bar y(s)) \in \partial\Omega$, one has $h\in C^1$ on a neighborhood of $(\bar y^0(s), \bar y(s))$ and there exists $\delta >0$ satisfying
\[
\begin{split}
  &\inf_{a \in A} \nabla_x h(\bar y^0(s), \bar y(s)) \cdot f(\bar y^0(s), \bar y(s), a) < -\delta , \\
  &\inf_{\w \in \C \cap \partial\B} \nabla_x h(\bar y^0(s), \bar y(s)) \cdot \sum_{i=1}^m g_i(\bar y^0(s), \bar y(s)) \w^i < -\delta;
\end{split}
\]
\item[]{\bf (IPFCn)}$_f$ for every $s\in [0,\bar S[$ such that $(\bar y^0(s), \bar y(s)) \in \partial\Omega$, one has $h\in C^1$ on a neighborhood of $(\bar y^0(s), \bar y(s))$ and there exists $\delta >0$ satisfying
\[
\begin{split}
  &\sup_{a \in A} \nabla_x h(\bar y^0(s), \bar y(s)) \cdot f(\bar y^0(s), \bar y(s), a) >\delta , \\
  &\sup_{\w \in \C \cap \partial\B} \nabla_x h(\bar y^0(s), \bar y(s)) \cdot \sum_{i=1}^m g_i(\bar y^0(s), \bar y(s)) \w^i >\delta,
\end{split}
\]
\end{itemize}}
respectively. Here $\nabla_x$ denotes the classical partial gradient operator w.r.t. the variable $x$.  
\end{remark}

 \begin{remark} \label{misuraE}
Note that for a feasible process  $\bar z=(\bar S, \bar\w^0 ,\bar\w,\bar\alpha , \bar y^0 ,\bar y , \bar\nu )$ with $h(\bar y^0(\bar s), \bar y(\bar s))=0$ for some $\bar s \in ]0,\bar S]$, we can assume that there exists some $\bar\varepsilon_1>0$ sufficiently small such that  $\ell(E(\bar s,\varepsilon))>0$ for  all $\varepsilon\in]0,\bar\varepsilon_1]$.
Indeed, if  $\exists\varepsilon >0$  such that $\ell(E(\bar s,\varepsilon)) =0$, then 
  $\ds \max_{(\xi_0,\xi) \in \partial^* h(\bar y^0(r), \bar y(r))} \left[ \xi_0  {\bar\w}^0(r) + \xi \cdot \bar F(r)   \right] < 0$   for a.e.  $r\in [\bar s-\varepsilon,\bar s].$   
But then, the function ${\mathcal H}:= h \circ (\bar y^0, \bar y)$ is differentiable a.e. in $[\bar s-\varepsilon, \bar s]$   and   verifies
\[ 
\begin{split}	
\frac{d{\mathcal H}}{ds}(s) &\leq  \max_{(\xi_0,\xi) \in \partial^C h(\bar y^0(s), \bar y(s))} \left[ \xi_0 \frac{d{\bar y}^0}{ds}(s) + \xi \cdot\frac{d{\bar y}}{ds}(s)   \right]   \\
&=\max_{(\xi_0,\xi) \in \partial^* h(\bar y^0(s), \bar y(s))} \left[ \xi_0 {\bar\w}^0(s) + \xi \cdot \bar F(s) \right]<0,
\end{split}
\]
since the scalar product is bilinear. Thus, for all $s \in [\bar s-\varepsilon, \bar s[$  one has
\[ -h(\bar y^0(s), \bar y(s))= h(\bar y^0(\bar s), \bar y(\bar s)) -h(\bar y^0(s), \bar y(s)) = \int_s^{\bar s}\frac{d{\mathcal H}}{ds}(\sigma)\,d\sigma <0,  \]
so that  $h(\bar y^0(s) , \bar y(s)) >0$, in contradiction with the feasibility of $\bar z$. In an analogous way,  one can derive that,  if $h(\bar y^0(\bar s), \bar y(\bar s))=0$ at   $\bar s \in [0,\bar S[$, then   $\ell(\Gamma(\bar s,\varepsilon))>0$ for  all $\varepsilon>0$ small enough.

\noindent Furthermore, by the bilinearity of the scalar product, in  \textbf{(CQn)$_b$},  \textbf{(CQn)$_f$},  and  \textbf{(CQn)$'_b$}, \textbf{(CQn)$'_f$},   one can replace  $\partial^*h(\bar y^0(s),\bar y(s))$ with  $\partial^C h(\bar y^0(s),\bar y(s))$. Hence, in particular,  all these conditions are verified   for any $(\xi_0,\xi)\in$ \linebreak $\partial^{>}h(\bar y^0(s),\bar y(s))$, since $\partial^{>}h(\bar y^0(s),\bar y(s))\subseteq \partial^C h(\bar y^0(s),\bar y(s))$.
\end{remark}

We can now establish the following normal Maximum Principle, which extends \cite[Th. 4.2]{MS20}.
\begin{theorem} \label{Th_CQN} Assume  {\bf (H0)}-{\bf (H2)}.
Let $(\bar S, \bar\w^0 ,\bar\w,\bar\alpha , \bar y^0 ,\bar y , \bar\nu )$ be an extended sense local minimizer for \eqref{Pe}. Assume that  $(p_0,p,\lambda,\pi,\mu)$,  $(m^0,m)$  meet  the conditions of Theorem \ref{Th_pmp} and verify the strenghtened non-triviality condition  \eqref{nondeg}. Then, 
\begin{itemize}
\item[{\rm (i)}] if  hypothesis {\rm\textbf{(CQn)$_b$}} is satisfied, one has
\begin{equation} \label{cqnb}
\begin{aligned}
&|q(\bar S)| + \lambda \neq 0 \qquad &\text{if } \bar y^0(\bar S) > \bar y^0(0), \\
&|q_0(\bar S)| + |q(\bar S)| + \lambda \neq 0 \qquad &\text{if } \bar y^0(\bar S) = \bar y^0(0).
\end{aligned}
\end{equation}
In particular, if {\rm proj}$_{(t_2,x_2)}N_{\T_0}(\bar y^0(0), \bar y(0), \bar y^0(\bar S), \bar y(\bar S)) = \{(\xi_{t_2},0) \}$ and $\xi_{t_2}=0$ whenever $\bar y^0(\bar S) =\bar y^0(0)$, then $\lambda \neq 0$;

\item[{\em (ii)}]  if  hypothesis {\rm\textbf{(CQn)$_f$}} is satisfied, one has
\begin{equation} \label{cqnf}
\begin{aligned}
&|q(0)| + \lambda \neq 0 \qquad &\text{if } \bar y^0(\bar S) > \bar y^0(0), \\
&|q_0(0)| + |q(0)| + \lambda \neq 0 \qquad &\text{if } \bar y^0(\bar S) = \bar y^0(0).
\end{aligned}
\end{equation}
 In particular, if {\rm proj}$_{(t_1,x_1)}N_{\T_0}(\bar y^0(0), \bar y(0), \bar y^0(\bar S), \bar y(\bar S)) = \{(\xi_{t_1},0) \}$ and $\xi_{t_1}=0$ whenever $\bar y^0(\bar S) =\bar y^0(0)$, then $\lambda \neq 0$.
 \end{itemize}
\end{theorem}
The proof of this result  is given in the Appendix.

\begin{remark}[Multiple state constraints] The normal Maximum Principle in Theorem \ref{Th_CQN} can be extended to the case with
 multiple state constraints $h_i$ for $i=1,\dots,N$, by modifying  condition {\bf (CQn)$'_b$} as follows: for every $s\in]0,\bar S]$ such that  $h_{i_j}(\bar y^0(s),\bar y(s))=0$ for  some indexes $1\leq i_1<\dots<i_k\leq N$,   $j=1,\dots,k$, there exist $\varepsilon$, $\delta>0$  and  a measurable control $(\tilde \w,\hat\alpha)$ taking values in $(\C\cap\partial \B)\times A$, satisfying    for all  $\sigma\in] s-\varepsilon, s[\cap[0,\bar S]$ and for all $j=1,\dots,k$:
\bel{cqn_multiple_constraints}
\ds \max_{(\xi_0,\xi)\in\partial^*h_{i_j}(\bar y^0(\sigma),\bar y(\sigma))} \xi\cdot \Big(F(\bar y^0,\bar y,\bar \w^0,\hat \w,\hat\alpha)(s') -\bar F(s')\Big)<-\delta,  
\eeq
where $\hat \w:=(1-\bar \w^0)\tilde \w$, for a.e. $s'\in \mathcal{E}(s,\varepsilon)$, defined as follows
\[
\begin{split}
&\mathcal{E}(s,\varepsilon):=\Bigl\{ r\in [s-\varepsilon,s]\cap[0, \bar S]\text{ : } \\
&\textcolor{white}{fdsfsdfsd}\max_{(\xi_0,\xi) \in \partial^*h_{i_j}(\bar y^0(r),\bar y(r))} [\xi_0 \bar\w^0(r)+\xi \cdot \bar F(r)]\ge0 \qquad \forall j=1,\dots,k    \Bigr\}.
\end{split}
\]
Notice that, arguing as in  Remark \ref{misuraE}, one can deduce that  $\ell(\mathcal{E}(s,\varepsilon))$ for all $\varepsilon>0$ sufficiently small.  Similar modifications are to be made for condition {\bf (CQn)$'_f$}.
\end{remark}
\vsm
 Theorem  \ref{Th_CQN} implies nondegenerate normality when essentially  the endpoint constraint either at the final or at the initial position is inactive. We provide below some sufficient conditions to guarantee normality even in some situations where initial and final positions lay on the boundary of the endpoint constraint.
\vsm
\noindent\textbf{Endpoint qualifications for normality} \textbf{(TQn)$_b$}, \textbf{(TQn)$_f$}. {\em  Let us consider  $\bar z:=(\bar S, \bar\w^0 ,\bar\w,\bar\alpha , \bar y^0 ,\bar y , \bar\nu )$, a feasible extended sense process for the extended optimization problem \eqref{Pe}.  
\begin{itemize}
\item[] {\bf (TQn)$_b$} We say that $\bar z$ meets condition  {\bf(TQn)$_b$} if there exists $\varepsilon>0$ such that $(\bar y^0(s), \bar y(s))\in\text{{\rm Int}}(\Omega)$ for each $s \in [\bar S-\varepsilon , \bar S[$ and one among the following conditions {\em (a)}, {\em (b)}  holds true:
\begin{equation} \label{inclusione1}
\begin{split}
\text{{\em (a)}} \quad &\left(-\text{{\rm proj}}_{(t_2,x_2)}N_{\T_0}(\bar y^0(0), \bar y(0),\bar y^0(\bar S), \bar y(\bar S)) \setminus \{ 0_{1+n}\} \right) \, \textcolor{white}{llllllllll} \\
&\textcolor{white}{fssdsdfsdfdf}\cap\,\partial^{>} h(\bar y^0(\bar S), \bar y(\bar S)) = \emptyset
\end{split}
\end{equation}
and for any $(\zeta_{t_2},\zeta_{x_2}) \in \text{{\rm proj}}_{(t_2,x_2)}N_{\T_0}(\bar y^0(0), \bar y(0),\bar y^0(\bar S), \bar y(\bar S)) $\linebreak $ +[0, +\infty[\cdot\partial^> h(\bar y^0(\bar S), \bar y(\bar S))$,\footnote{As in Remark \ref{RMSC2}, we adopt the convention that $X+\emptyset=X$.}  
  one has
\begin{equation} \label{tqn1}
\min_{a \in A} \left[\zeta_{x_2} \cdot f(\bar y^0(\bar S), \bar y(\bar S),a) +\zeta_{t_2}  \right] <0 \qquad \text{if } (\zeta_{t_2},\zeta_{x_2}) \neq (0,0);
\end{equation}
\begin{equation} \label{inclusione2}
\begin{split}
\text{{\em (b)}}\quad &\Big(-\text{{\rm proj}}_{x_2}N_{\T_0}(\bar y^0(0), \bar y(0),\bar y^0(\bar S), \bar y(\bar S)) \setminus \{ 0_{n}\} \Big) \, \textcolor{white}{lllllllllllllllll}\\
&\textcolor{white}{fssdsdfs} \cap \,\text{{\rm proj}}_{x}\partial^> h(\bar y^0(\bar S), \bar y(\bar S))= \emptyset 
\end{split}
\end{equation}
and for any $(\zeta_{t_2},\zeta_{x_2}) \in \text{{\rm proj}}_{(t_2,x_2)}N_{\T_0}(\bar y^0(0), \bar y(0),\bar y^0(\bar S), \bar y(\bar S))$\linebreak $+[0, +\infty[\cdot \partial^> h(\bar y^0(\bar S), \bar y(\bar S))$ with $\zeta_{x_2}\neq0$, one has
\begin{equation} \label{tqn2}
\bar y^0(\bar S) > \bar y^0(0), \,\,\, \bar\nu(\bar S) <K, \,\,\, \min_{\omega \in \C \cap \partial \B} \left[ \zeta_{x_2} \cdot \left( \sum_{j =1}^m g_j(\bar y^0(\bar S), \bar y(\bar S)) \omega^j  \right)   \right] <0.
\end{equation}

\item[]{\bf (TQn)$_f$}  We say that $\bar z$ meets condition  {\bf(TQn)$_f$} if there exists $\varepsilon>0$ such that $(\bar y^0(s), \bar y(s))\in\text{{\rm Int}}(\Omega)$ for each $s \in ]0,\varepsilon]$ and one among the following conditions {\em (a)}, {\em (b)}  holds true:
\begin{equation} \label{inclusione3}
\begin{split}
\text{{\em (a)}} \quad &\left(-\text{{\rm proj}}_{(t_1,x_1)}N_{\T_0}(\bar y^0(0), \bar y(0),\bar y^0(\bar S), \bar y(\bar S)) \setminus \{ 0_{1+n}\} \right) \, \textcolor{white}{llllllllll} \\
&\textcolor{white}{fssdsdfsdfdf}\cap\,\partial^{>} h(\bar y^0(0), \bar y(0)) = \emptyset
\end{split}
\end{equation}
and  for any $(\zeta_{t_1},\zeta_{x_1}) \in \text{{\rm proj}}_{(t_1,x_1)}N_{\T_0}(\bar y^0(0), \bar y(0),\bar y^0(\bar S), \bar y(\bar S))$ \linebreak $+[0, +\infty[\cdot\partial^> h(\bar y^0(0), \bar y(0))$ one has
\begin{equation} \label{tqn3}
\max_{a \in A} \left[ \zeta_{x_1} \cdot f(\bar y^0(0), \bar y(0),a) +\xi_{t_1}  \right] >0 \qquad \text{if }  (\zeta_{t_1},\zeta_{x_1})  \neq (0,0);
\end{equation}
\begin{equation} \label{inclusione4}
\begin{split}
\text{{\em (b)}}\quad &\Big(-\text{{\rm proj}}_{x_1}N_{\T_0}(\bar y^0(0), \bar y(0),\bar y^0(\bar S), \bar y(\bar S)) \setminus \{ 0_{n}\} \Big) \,  \textcolor{white}{lllllllllllllllll}\\
&\textcolor{white}{fssdsdfs} \cap \,\text{{\rm proj}}_{x}\partial^> h(\bar y^0(0), \bar y(0))= \emptyset 
\end{split}
\end{equation}
and   for any $(\zeta_{t_1},\zeta_{x_1}) \in \text{{\rm proj}}_{(t_1,x_1)}N_{\T_0}(\bar y^0(0), \bar y(0),\bar y^0(\bar S), \bar y(\bar S))$ \linebreak $+[0, +\infty[\cdot\partial^> h(\bar y^0(0), \bar y(0))$ with $\zeta_{x_1}\neq0$, one has
\begin{equation} \label{tqn4}
\bar y^0(\bar S) > \bar y^0(0),\ \  \bar\nu(\bar S) <K, \ \  \max_{\omega \in \C \cap \partial \B}   \zeta_{x_1} \cdot   \sum_{j =1}^m g_j(\bar y^0(0), \bar y(0)) \omega^j    >0.
\end{equation}
\end{itemize}}
Condition  {\bf (TQn)$_b$} generalizes the endpoint constraint qualifications considered in  \cite{MS20} for the case with fixed initial point, which were in turn inspired by no gap conditions in \cite{MRV,AMR15}.  Notice that both  conditions 
 \eqref{inclusione1}, \eqref{inclusione2} [resp.,   \eqref{inclusione3},  \eqref{inclusione4}]
 are trivially verified whenever $(\bar y^0(\bar S),\bar y(\bar S)) \in\text{{\rm Int}}(\Omega)$ [resp., $(\bar y^0(0),\bar y(0))\in\text{{\rm Int}}(\Omega)$], since  $\partial^{>} h(\bar y^0(\bar S), \bar y(\bar S)) =\emptyset$ [resp., $ \partial^{>} h(\bar y^0(0), \bar y(0))=\emptyset $]. 
 
\begin{proposition} \label{TQN}  Assume  {\bf (H0)}-{\bf (H2)}.
Let $(\bar S, \bar\w^0 ,\bar\w,\bar\alpha , \bar y^0 ,\bar y , \bar\nu )$ be an extended sense local minimizer for \eqref{Pe}. Assume that  $(p_0,p,\lambda,\pi,\mu)$,  $(m_0,m)$  meet  the conditions of Theorem \ref{Th_pmp}. Then,  when either {\em (i)} or {\em (ii)} below holds true, one has $\lambda\ne0$:
\begin{itemize} 
\item[{\rm(i)}] hypothesis {\bf(TQn)$_b$} is satisfied and the multipliers $(p_0,p,\lambda,\pi,\mu)$ verify the strengthened non-triviality condition \eqref{cqnb}; 
\item[{\rm (ii)}] hypothesis {\bf(TQn)$_f$} is satisfied and the multipliers $(p_0,p,\lambda,\pi,\mu)$ verify the strengthened non-triviality condition \eqref{cqnf}.
\end{itemize}
\end{proposition}
The proof of this result  is postponed to the Appendix.
\begin{remark}[Multiple state constraints] Proposition \ref{TQN} can be easily adapted to the case with multiple state constraints  $h_i$ for $i=1,\dots,N$ by simply choosing $h:=  h_1\vee\dots\vee h_N$. In particular,   by applying the max-rule for subdifferentials, in \textbf{(TQn)}$_b$, (b),  condition \eqref{inclusione2} can be replaced with 
\[
\begin{split}
&\left(\sum_{i=1}^N\, \text{{\rm proj}}_{x}[0,+\infty[\cdot\partial^> h_i(\bar y^0(\bar S), \bar y(\bar S))\right) \,\\
&\textcolor{white}{dfsfsd} \cap\, \left(-\text{{\rm proj}}_{x_2}N_{\T_0}(\bar y^0(0), \bar y(0),\bar y^0(\bar S), \bar y(\bar S)) \setminus \{ 0_{n}\} \right) = \emptyset,
\end{split}
\]
and one can require that  condition  (\ref{tqn2}) is satisfied  
for all $(\zeta_t,\zeta_x)$ in the set $\text{{\rm proj}}_{(t_1,x_1)}N_{\T_0}(\bar y^0(0), \bar y(0),\bar y^0(\bar S), \bar y(\bar S))$   $+ \sum_{i=1}^N \, [0, +\infty[\cdot\partial^> h_i(\bar y^0(\bar S), \bar y(\bar S)). $
The  other conditions can be adapted in a similar way.
\end{remark}

\vsm
From Propositions  \ref{nondeg_norm},  \ref{TQN}, and Theorem \ref{Th_CQN}  we  deduce as a corollary the main result of this section.  

\begin{theorem}\label{Th_gapsuff} Assume  {\bf (H0)}-{\bf (H2)}.
 Consider the optimal control problem \eqref{P} and its extended sense formulation \eqref{Pe}.   Assume that there exists a local or global extended sense minimizer $(\bar S, \bar\w^0 ,\bar\w,\bar\alpha , \bar y^0 ,\bar y , \bar\nu )$ such that condition  {\rm \textbf{(CNa)}} and either    {\rm\textbf{(CQn)$_b$}}-{\rm\textbf{(TQn)$_b$}} or  {\rm\textbf{(CQn)$_f$}}-{\rm\textbf{(TQn)$_f$}}  are verified.   
Then $(\bar S, \bar\w^0 ,\bar\w,\bar\alpha , \bar y^0 ,\bar y , \bar\nu )$  is a  normal extremal and, in consequence of Theorem \ref{Th_Norm},  there is no local or global  infimum gap, respectively.
\end{theorem}

\section{Some examples}\label{S4}
Let us illustrate the preceding theory through some  examples.

\begin{example}\label{es1}   In this example   the absence of an infimum gap can be easily deduced from the sufficient  conditions introduced in Section \ref{S3}.       Consider the problem
\bel{strictex3}
\left\{
\begin{array}{l}
\mbox{minimize }\  -x(1)
\\ [1.5ex]
\mbox{over } (x,v, u)=(x^1,x^2,x^3,v, u^1,u^2) \in W^{1,1}([0,1];\R^{3}\times\R\times\R^2)   \quad
\mbox{satisfying }\\ [1.5ex]
\begin{array}{l}
\ds \frac{d x}{dt}(t) = f(x(t))+g_1(x(t))\,\frac{d u^1}{dt}(t) +g_2(x(t))\,\frac{d u^2}{dt}(t) 
\\ [1.5ex]
\ds \frac{d v}{dt}(t) =  \, \left|\frac{d u}{dt}(t)\right|  \\ [1.5ex]
\displaystyle\frac{du}{dt}(t)\in\C:=\R^2  \qquad  \mbox{ a.e. } t \in [0,1],   \\ [1.5ex]
x(t)\in\Omega(t)  \qquad\forall t\in[0,1],
 \\ [1.5ex]
\displaystyle v(0)=0, \ \ v(1)   \leq 2 ,\quad 
 x(0)=\T^1_0, \quad   x(1)\in\T^2_0 .\end{array}
\end{array}
\right.
\eeq
in which 
$$
\begin{array}{l}
 \Omega(t):=\{(t,x)\in\R^{4}: \ -1\le x^1\le1+t, \  -1\le x^2 \le 1, \   -1\le x^3 \le 1  \}, \\[1.5ex]
  \ \T^1_0:= \{x\in\R^3: \   (x^1-1)^2+(x^2)^2+(x^3)^2\le 1/9, \ x^1\le 1  \},   \\[1.5ex]
  \ \T^2_0:=  \{x\in\R^3: \   (x^1+1)^2+(x^2)^2+(x^3)^2\le 1, \ x^1\ge -1  \},  
\end{array}
$$
 and
$$
\ds g_1(x):=\left(\begin{array}{l}  1 \\   0\\0\end{array}\right),   \quad \ds g_2(x):=\left(\begin{array}{l} \ \ 0 \\ -1\\-x^1\end{array}\right) , \quad  \ds f(x):=\left(\begin{array}{l} \ 0 \\  x^2x^3\\ \ 0\end{array}\right)  \quad \forall  x\in \R^3\,.
$$
The extended  problem is
\bel{extx3}
\left\{
\begin{array}{l}
\mbox{minimize } \  -y^1(S)
\\ [1.5ex]
\mbox{over } S>0, \ (y^0, y^1,y^2,y^3,\nu, \w^0,\w^1,\w^2) \in W^{1,1}([0,S])   \quad
\mbox{satisfying }\\ [1.5ex]
\begin{array}{l}
\ds \frac{d y^0}{ds}(s) = \w^0(s)
\\ [1.5ex]
\ds \frac{d y }{ds}(s) =f(y(s))\w^0(s)+g_1(y(s))\,\w^1(s) +g_2(y(s))\,\w^2(s)
\\ [1.5ex]
\ds \frac{d \nu}{ds}(s) =  \, \left|\w(s)\right|  \\ [1.5ex]
\displaystyle (\w^0,\w)(s) \in\CC \qquad  \mbox{ a.e. } s \in [0,S],   \\ [1.5ex]
y(s)\in\Omega(y^0(s))  \qquad\forall s\in[0,S],
 \\ [1.5ex]
\ds \nu(0)=0, \,\, \nu(S)   \leq 2 ,\,\, 
 (y^0(0), y(0))=\{0\}\times \T^1_0, \,\,  (y^0(S), y(S))\in\{1\}\times \T^2_0 .\end{array}
\end{array}
\right.
\eeq
An extended sense  minimizer  is clearly given by the feasible extended sense process   $(\bar S, \bar\w^0 ,\bar\w,  \bar y^0 ,\bar y , \bar\nu )$, where  
\bel{optc}
\bar S=2, \qquad (\bar\w^0,\bar\w)=(\bar\w^0,\bar\w^1,\bar\w^2)  =(1,0,0)\chi_{_{[0,1]}}  +(0,-1,0)\chi_{_{]1,2]}}   \,,
\eeq
  and  one considers the corresponding trajectory with initial state condition $y(0)=(1,0,0)$, namely, 
\bel{optt}
(\bar y^0,\bar y,\bar\nu) =(\bar y^0,\bar y^1,\bar y^2,\bar y^3,\bar\nu)  =(s,1,0,0,0 )\chi_{[0,1]} +(1,2-s,0,0,s-1)\chi_{[1,2]} \,.
\eeq
 It is not difficult  to check that this process verifies   conditions  {\bf (CNa)}, {\bf (CQn)$_b$},  and {\bf (TQn)$_b$}. In consequence, the absence of a  gap between the infima of  problems \eqref{strictex3}, \eqref{extx3} follows  directly from Theorem \ref{Th_gapsuff}.
   \end{example}
 
Next example shows how the criterion of normality can guarantee the absence of the infimum gap in situations where other sufficient conditions fail.

\begin{example}\label{es2}  Consider again the minimization problem \eqref{strictex3} and its extended version \eqref{extx3}, where $\T^1_0$ is as above, while the time-dependent state constraint $\Omega(t)$ and the final-point constraint $\T^2_0$ are replaced with 
$$
\begin{array}{l}
\Omega:= \{(x^1,x^2,x^3): \  -1\le x^1  \le 1, \ -1\le x^2  \le 1, \ -1\le x^3  \le 1 \}, \\[1.5ex]
 \T^2_0:=  \{(x^1,x^2,x^3): \ -1\le x^1\le 0, \ 0\le x^2\le1, \ 0\le x^3\le 1\},
 \end{array}
 $$
 respectively.   Then the extended sense process  $(\bar S, \bar\w^0 ,\bar\w,  \bar y^0 ,\bar y , \bar\nu )$ given by \eqref{optc}, \eqref{optt} is still a (feasible) minimizer for the extended problem \eqref{extx3}. However, as it is easy to check, now the nondegeneracy condition  {\bf (CNa)} is met,  but none of the conditions  {\bf (CQn)$_b$}, {\bf (TQn)$_b$}, {\bf (CQn)$_f$}, and {\bf (TQn)$_f$} is satisfied.

\noindent From  Theorem  \ref{Th_pmp} there exist a set of multipliers $(p_0 ,p,\pi, \lambda,\mu)$ and functions $(m_0,m)$  with  $\pi=0$,  since $\nabla_v\Psi\equiv0$ and  $\bar\nu(2)=1<2$. Also, $m_0\equiv 0$ as the state constraint does not depend on time, $\mu([0,2])=\mu([0,1])$, and $m(s)\in\partial^>_xh(\bar y(0))$ $\mu$-a.e.  yields  $m(s)=(1,0,0)$ $\mu$-a.e. in $[0,1]$.  By the adjoint equation it follows that  the path $(p_0,p)=(p_0,p_1,p_2,p_3)\equiv(\bar p_0,\bar p_1,\bar p_2,\bar p_3)$ is constant. 
From  the transversality conditions 
\bel{ex3f2}
\begin{split}
&(p_0,p_1,p_2,p_3)(0)\in \R\times N_{\T^1_0}(1,0,0), \\ 
&-(q_0,q_1,q_2,q_3)(2)\in\lambda\{(0,-1,0,0)\}+\R\times N_{\T^2_0}(0,0,0), 
\end{split}
\eeq
where $q_0\equiv\bar p_0$,   and $q(s)=(\bar p_1+\mu([0,1]), \bar p_2,\bar p_3)$ for all $s\in]1,2]$, we derive that $\bar p_0\in\R$, $\bar p_1\ge0$, $\bar p_2=\bar p_3=0$, $q_1(2)=\bar p_1+\mu([0,1])=\lambda -\alpha_1 $ with $\alpha_1\ge0$.  The maximality condition implies the relations 
\bel{ex3f3}
\begin{array}{l}
 \bar p_0\chi_{[0,1]}(s) =0, \quad -q_1(s)\chi_{]1,2]}(s)=0,
 \end{array}
\eeq
from which  we deduce that $\bar p_0=0$ and $\bar p_1+\mu([0,1])=\lambda -\alpha_1=0 $.   Hence, recalling that $\bar p_1\ge0$, we get $\bar p_1=\mu([0,1])=0$, $\lambda=\alpha_1\ge0$. 
So, the strengthened non-triviality  condition  
$ \|p\|_{L^\infty}+\mu([0,2])+\lambda\ne0$
 implies that $\lambda\ne0$ and this shows that  $(\bar S, \bar\w^0 ,\bar\w,  \bar y^0 ,\bar y , \bar\nu )$ is a normal extremal.  Consequently, there is no infimum gap for the `normality test' established in Theorem  \ref{Th_Norm}.
  \end{example}

 However,  normality itself is only a sufficient condition to avoid  the gap  (even for systems with   drift, $f$).
 
\begin{example}\label{es3}  
Let us consider the problem in Example \ref{es2} where we only modify the initial-point  target $\T^1_0$, replacing it with   $\T^1_0:=\{(1,0,0)\}$. Then the extended sense process   $(\bar S, \bar\w^0 ,\bar\w,  \bar y^0 ,\bar y , \bar\nu )$ of before is obviously still admissible and minimizing, but it is easy to see that the set of degenerate multipliers $(p_0,p,\lambda,\mu)$ with $p_0=p_2=p_3=0$, $p_1=-1$,   $\mu=\delta_{\{0\}}$, and $\lambda=0$   meets all the conditions of   Theorem   \ref{Th_pmp}. So,   $(\bar S, \bar\w^0 ,\bar\w,  \bar y^0 ,\bar y , \bar\nu )$ is an abnormal extremal. But there is no gap, because, for any $n\in\N$, the strict sense process   $(\bar S, \w^0_n ,\w_n,   y^0_n , y_n ,\nu_n )$, where 
$$
\bar S=2, \qquad (\w^0_n,\w_n)=    \left(1-\frac{1}{n},-\frac{1}{n},0\right)\chi_{_{[0,1]}}  +\left(\frac{1}{n},-1+\frac{1}{n},0\right)\chi_{_{]1,2]}}   \,,
$$
  and  $(y^0_n , y_n ,\nu_n )$ is the corresponding trajectory of the control system in \eqref{extx3} with initial   condition $(y^0,y,\nu)(0)=(0,1,0,0,0)$, is feasible and minimizing for the original problem, since it has cost  equal to zero. 
\end{example}

\appendix
\section{ Appendix}
\noindent {\it Proof of Theorem \ref{Th_CQN}: }\,  Thanks to Remark \ref{misuraE}, when  \textbf{(CQn)$_b$} is in force the proof of Theorem \eqref{Th_CQN} is analogous to the proof of \cite[Th. 4.2]{MS20},  while under assumption   \textbf{(CQn)$_f$} it  requires some adaptation. For this reason, we limit ourselves to  give the proof in the last case. 

\noindent By standard truncation and mollification arguments, we can assume   $h$  Lipschitz continuous, with Lipschitz constant $L>0$, and   $f$,  $g_1,\dots,g_m$, and their  limiting subdifferentials in $(t,x)$,  $L^\infty$-bounded  by some constant  $\tilde M >0$.  Set $M:=(1+m)\tilde M$.

\noindent By assumption, the local minimizer $(\bar S, \bar\w^0 ,\bar\w,\bar\alpha , \bar y^0 ,\bar y , \bar\nu )$ has a set of multipliers $(p_0,p,\pi,\lambda, \mu)$ and some functions $(m_0,m)$ such that the conditions (i)-(vi) of Theorem \ref{Th_pmp} hold true, and verifying the strengthened non-triviality condition (\ref{nondeg}).
Let us first assume that $\bar y^0(\bar S)>\bar y^0(0)$ and suppose by contradiction that
\begin{equation} \label{hpcontradd}
q(0)=0, \qquad \lambda=0.
\end{equation}
Set
\[
\bar s:=\sup\{s\in[0,\bar S]: \ \ \mu(]0,s])=0\}.
\]
Observe that $\bar s<\bar S$.  Indeed, if not, $\mu(]0,\bar S])=0$. But in this case $q(s)=p(s)+\mu(\{0\})m(0)$, so that it is absolutely continuous and by the adjoint equation  with initial condition $q(0)=0$ it follows that  $q\equiv 0$. 
Precisely, by  known properties of the convex hull of the limiting subdifferential of locally Lipschitz continuous functions (see e.g. \cite[Ch. 6]{OptV}),  we have
\[	|q(s)|\le \int_0^s\left|  \frac{dq}{ds} (s)\right|\,ds   \leq M  \int_0^s |q(s)|\, ds,		\]
 which implies that $q\equiv 0$ by Gronwall's Lemma. Since $\lambda=0$ by \eqref{hpcontradd}, this is in  contradiction with  the first relation in \eqref{nondeg}. 
 When $\bar y^0(\bar S)=\bar y^0(0)$ and we assume  by contradiction that
\begin{equation} \label{hpcontradd1}
q_0(0)=0, \qquad q(0)=0, \qquad \lambda=0,
\end{equation}
the value $\bar s$ defined as above is still strictly smaller than $\bar S$, since otherwise $\mu(]0,\bar S])=0$,  so that  $(q_0,q)\equiv 0$, again  by the adjoint equation.  In view of \eqref{hpcontradd1}, this yields contradiction with the second relation in \eqref{nondeg}.
Obviously, $(\bar y^0(\bar s),\bar y(\bar s))\in\partial\Omega$. 

\vsm From now on, the proof is the same for both cases. 
Introduce 
$$(z_0,z)(s):=\left(p_0(s)+m_0(0)\mu(\{0\}), p(s)+m(0)\mu(\{0\})\right),
$$
  so that,  for any $s\in[0,\bar S[$,
\begin{equation} \label{intermedio1}
\begin{array}{l}
\ds(q_0(s),q (s))=\left(p_0 (s)+\int_{[0,s[}m_0(\sigma)\mu(d\sigma)\,, \, p(s)+\int_{[0,s[}m(\sigma)\mu(d\sigma)\right) \\
\ds\qquad\qquad\qquad\qquad\qquad =\left(z_0(s)+\int_{]0,s[}m_0(\sigma)\mu(d\sigma)\,,\, z(s)+\int_{]0,s[}m(\sigma)\mu(d\sigma)\right).
\end{array}
\end{equation}
By the adjoint equation, $(z_0,z)$ verifies
\begin{equation} \label{intermedio2}
\begin{cases} 
\ds-\left(\frac{dz_0}{ds}, \frac{dz}{ds}\right)(s)
 \in \text{co}\,\partial_{t,x}\{q(s) \cdot \bar F(s)  \} = \text{co}\, \partial_{t,x} \left\{z(s)\cdot \bar F(s) + \int_{]0,s[} m(\sigma) \mu(d\sigma) \cdot \bar F(s) \right\} \\
z(0)=0, \quad  (\text{and $z_0(0)=0$, \, if  $\bar y^0(\bar S)=\bar y^0(0)$}).
\end{cases}
\end{equation}
Since the integral on the right hand side is identically zero in $]0,\bar s[$, arguing as above  we derive that $z(s)=0$ and therefore $q(s)=0$ for all $s\in[0,\bar s[$,  by continuity. Moreover,   Gronwall's Lemma implies that $|z(s)|\leq C\,\mu([\bar s,s[)$ for all $s\in[\bar s,\bar S[$, for some $C>0$, so   that 
\begin{equation} \label{intermedio3}
|q(s)|\leq |z(s)|+L \mu([\bar s,s[)\leq ( C+L) \mu([\bar s,s[) \qquad \forall s\in[\bar s,\bar S[.
\end{equation}
 As a consequence of (\ref{intermedio1}), for every $s\in[\bar s,\bar S]$   one gets $q(s) =z(s) + \int_{[\bar s,s[}m(\sigma) \mu(d\sigma)$,  and (\ref{intermedio2}), (\ref{intermedio3}) imply
\begin{equation} \label{valfinq}
\left|q(s)-\int_{[\bar s,s[}m(\sigma)\mu(d\sigma)\right| \leq  M \int_{\bar s}^s |q(\sigma)| \,d\sigma      \leq   \bar C\,\mu([\bar s,s[)(s-\bar s),
\end{equation}
where $\bar C:=M(L+C)$. 
In view of  {\bf(CQn)$'_f$} in Remark \ref{Rem2cqn},  there exist $\varepsilon$, $\delta>0$ and a  measurable control $(\tilde \w,\hat \alpha): [0,\bar S]\to (\C\cap\partial\B)\times A$,  verifying for all  $(\xi_0,\xi)\in\partial^*h(\bar y^0(s),\bar y(s))$  with $s\in]\bar s,\bar s+\varepsilon[\cap[0,\bar S]$:
\begin{equation} \label{centr}
  \xi\cdot \left(F((\bar y^0,\bar y,\hat \w^0,\hat \w, \hat\alpha)(s'))-\bar F(s')\right)> \delta, \quad\text{for a.e. $s'\in \Gamma(\bar s,\varepsilon)$,}
\end{equation}
where
$
(\hat \w^0(s),\hat \w(s),\hat\alpha(s)):=\left(\bar \w^0(s), (1-\bar \w^0(s))\tilde \w(s), \hat\alpha(s)\right)$  for a.e. $s\in [0,\bar S]$.  
Observe that, being  $\hat \w^0\equiv \bar \w^0$, one has   $|\hat \w|=1-\bar \w^0=|\bar \w|$ a.e.
As observed in  Remark \ref{misuraE},  $\ell(\Gamma(\bar s,\varepsilon))$ is $>0$ for any $\varepsilon >0$ sufficiently small. Moreover, \eqref{centr} is valid  for  any   $(\xi_0,\xi)\in \partial^C h(\bar y^0(s),\bar y(s))$, so that it is true, in particular,    for $(\xi_0,\xi)\in \partial^{>}h(\bar y^0(s),\bar y(s))$.
 On the other hand, by the maximization condition (\ref{maxham})  of Theorem \ref{Th_pmp}, it follows that,  for a.e. $s\in ]\bar s,\bar s+\varepsilon[\cap[0,\bar S]$,
\begin{equation} \label{fe3dim}
\begin{split}
&q_{0}(s)\left(\hat \w^0(s)-\bar \w^0(s)\right)+ \pi(|\hat \w(s)|-|\bar \w(s)|)  
+q(s)\left[F((\bar y^0,\bar y,\hat \w^0, \hat \w, \hat \alpha)(s))-\bar F(s)\right]  \\
&\textcolor{white}{fdsjhfsdkhksdj}=q(s)\left[F((\bar y^0,\bar y,\hat \w^0,\hat \w,\hat \alpha)(s))-\bar F(s)\right]\leq 0.
\end{split}
\end{equation}
Putting together  (\ref{valfinq}), (\ref{centr}), and  (\ref{fe3dim})  we get the desired contradiction. Indeed, for $\varepsilon>0$ small enough,  for any  $s'\in \Gamma(\bar s, \varepsilon)$,   one has 
\[
\begin{split}
0 &\geq q(s')\left[F((\bar y^0,\bar y,\hat \w^0,\hat \w,\hat \alpha)(s'))-\bar F(s')\right]  \\
&=\left( q(s')-\int_{[\bar s,s'[}m(\sigma)\mu(d\sigma)+ \int_{[\bar s,s'[}m(\sigma)\mu(d\sigma)\right) \left[F((\bar y^0,\bar y,\hat \w^0,\hat \w,\hat \alpha)(s'))-\bar F(s')\right]  \\
&\geq  \int_{[\bar s,s'[} m(\sigma)\,\left[F((\bar y^0,\bar y,\hat \w^0,\hat \w,\hat \alpha)(s'))-\bar F(s')\right]  \mu(d\sigma)-2M\bar C\,\mu([\bar s,s'[)(s'-\bar s) \\
&\geq \mu([\bar s,s'[) \left[ \delta -2M\, \bar C\,(s'-\bar s) \right]> 0
\end{split}
\]
for $\varepsilon>0$ sufficiently small. This concludes the proof.
\qed
\vsm

\noindent {\it Proof Proposition \ref{TQN}: }\,
The proof follows the same lines of the proof of  \cite[Prop. 4.1]{MS20}, where however only condition \textbf{(TQn)$_b$} for an implicit  state constraint $\Omega\subseteq\R^n$  is considered.  
\vsm
\noindent Let us prove (i).  Assume by contradiction  $\lambda=0$. Then
  the  transversality condition  \eqref{trans_cond} implies that
$$
( -q_0(\bar S), -q(\bar S))=(\zeta_{t_2},\zeta_{x_2}) \in \text{{\rm proj}}_{(t_2,x_2)}N_{\T_0}(\bar y^0(0), \bar y(0),\bar y^0(\bar S), \bar y(\bar S)), 
$$
where  $(\zeta_{t_2},\zeta_{x_2})\ne(0,0)$ and, in particular,  $\zeta_{x_2}\ne0$ if $\bar y^0(\bar S)>\bar y^0(0)$ by \eqref{cqnb}.   By hypothesis \textbf{(TQn)$_b$},  there is some $\varepsilon>0$ such that $(\bar y^0(s),\bar y(s))\in\,$Int$( \Omega)$ for all $s \in [\bar S-\varepsilon, \bar S[$. Hence $\mu([\bar S-\varepsilon,\bar S[)=0$, so that,  for any $s\in]\bar S-\varepsilon,\bar S[$,    $(q_0,q)$ is continuous at $s$ and
$$
(q_0(s),q(s))=\Big(p_0(s)+\int_{[0,\bar S-\varepsilon]} m_0(r)\mu(dr)\,,\,p(s)+\int_{[0,\bar S-\varepsilon]} m(r)\mu(dr)\Big).
$$ 
 Set 
$
(q_0(\bar S^-),q(\bar S^-)):=\lim_{s\to\bar S^-}(q_0(s),q(s))=(p_0,p)(\bar S)+\int_{[0,\bar S-\varepsilon]} (m_0,m)(r)\mu(dr). 
$
We get
$$
\begin{array}{l}
(q_0(\bar S^-),q(\bar S^-))=\Big(q_0(\bar S)-m_0(\bar S)\mu(\{\bar S\})\,,\,q(\bar S)-m(\bar S)\mu(\{\bar S\})\Big)
 = ( -\tilde\zeta_{t_2}, -\tilde\zeta_{x_2}),
\end{array}
$$
where  $(\tilde\zeta_{t_2},\tilde\zeta_{x_2}):=\Big(\zeta_{t_2}+\mu(\{\bar S\})m_0(\bar S),\zeta_{x_2}+\mu(\{\bar S\})m(\bar S)\Big)$. Thus, in particular,   the pair $(\tilde\zeta_{t_2},\tilde\zeta_{x_2})$ verifies
\begin{equation}\label{ip2}
(\tilde\zeta_{t_2},\tilde\zeta_{x_2})\in \text{{\rm proj}}_{(t_2,x_2)}N_{\T_0}(\bar y^0(0), \bar y(0),\bar y^0(\bar S), \bar y(\bar S))+ [0, +\infty[\cdot\partial^> h(\bar y^0(\bar S), \bar y(\bar S)).
\end{equation}
The continuity of $(q_0,q)$ on  $]\bar S-\varepsilon,\bar S[$ also implies that the equality (\ref{maxham}) in the Maximum Principle is verified for all  $s\in]\bar S-\varepsilon,\bar S[$. Hence, passing to the limit in it as $s$ tends to $\bar S^-$, we obtain
\begin{equation} \label{f31}
\left(\tilde\zeta_{x_2} \cdot f(\bar y^0(\bar S), \bar y(\bar S), a) + \tilde\zeta_{t_2}  \right) \w^0 + \tilde\zeta_{x_2} \cdot  \sum_{j=i}^m g_j(\bar y^0(\bar S), \bar y(\bar S)) \w^j   -\pi |\w| \geq 0.
\end{equation}
Suppose first that condition (a) in \textbf{(TQn)$_b$} is satisfied.  Then, from  \eqref{inclusione1}  we deduce  that   $(\tilde\zeta_{t_2},\tilde\zeta_{x_2}) \neq (0,0)$ and  choosing  $\w=0$ in \eqref{f31} we obtain a contradiction to \eqref{tqn1}.  

\noindent  If instead  condition (b) in \textbf{(TQn)$_b$} is valid,  $\pi=0$ and \eqref{cqnb}  implies that   $\zeta_{x_2}\ne0$. In view of  \eqref{ip2} and  hypothesis  \eqref{inclusione2},  this yields   $\tilde\zeta_{x_2} \neq 0$.  At this point,  we get a contradiction  to \eqref{tqn2} by choosing $\w^0=0$ in \eqref{f31}. 
 
The proof of  (ii) is very similar, hence we omit it.
\qed
\vsm


%
%



\end{document}